\documentclass[11pt]{article}

\usepackage{amsmath,amsfonts,amssymb,amsthm,mathtools,bbm, geometry}

\newtheorem{thm}{Theorem}[section]

\newtheorem{dfn}[thm]{Definition}

\newtheorem{lem}[thm]{Lemma}

\newtheorem{prop}[thm]{Proposition}

\newtheorem{remark}[thm]{Remark}

\newtheorem{cor}[thm]{Corollary}

\newtheorem{ex}[thm]{Example}

\newtheorem{question}[thm]{Question}

\def\bq{\begin{question}}
\def\bt{\begin{thm}}
\def\bp{\begin{prop}}
\def\blem{\begin{lem}}
\def\bd{\begin{dfn}}
\def\br{\begin{remark}}
\def\bc{\begin{cor}}
\def\bex{\begin{ex}}
\def\beqs{\begin{eqnarray*}}
\def\beq{\begin{eqnarray}}
\def\bi{\begin{itemize}}

\def\eq{\end{question}}
\def\et{\end{thm}}
\def\ep{\end{prop}}
\def\elem{\end{lem}}
\def\ed{\end{dfn}}
\def\er{\end{remark}}
\def\ec{\end{cor}}
\def\eex{\end{ex}}
\def\eeqs{\end{eqnarray*}}
\def\eeq{\end{eqnarray}}
\def\ei{\end{itemize}}

\def\ds{\displaystyle}

\def\c{\cdot}

\def\ov{\overline}
\def\r{\rangle}
\def\l{\langle}

\def\H{{\cal H}}
\def\K{{\cal K}}

\def\Hp{{\cal H}_\pi}

\def\w*{w^*-w^*}

\def\F{{\cal F}}

\def\ra{\rightarrow}

\def\ja{j_\alpha}

\def\xpe{\xi*_\pi \eta}

\def\vp{\varphi}
\def\DBG{\Delta(B(G))} 
\def\bs{\backslash} 

\def\bB{{\mathbb B}}
\def\bD{\Delta\!\!\!\!\Delta}

\def\O{\mathbb{O}}
\def\Q{\mathbb{Q}}

\def\gap{G^{ap}}
\def\dap{\delta^{ap}} 
\def\tdap{\widetilde{\delta^{ap}}}
\def\eap{e_{ap}}

    \def\AFG{A_\mathbbmtt{F}(G)}
    
       \def\pif{{\pi_\mathbbmtt{F}}} 

\def\PDA{{\cal P}_{\Delta(A)}}

\makeatletter
\def\moverlay{\mathpalette\mov@rlay}
\def\mov@rlay#1#2{\leavevmode\vtop{%
   \baselineskip\z@skip \lineskiplimit-\maxdimen
   \ialign{\hfil$\m@th#1##$\hfil\cr#2\crcr}}}
\newcommand{\charfusion}[3][\mathord]{
    #1{\ifx#1\mathop\vphantom{#2}\fi
        \mathpalette\mov@rlay{#2\cr#3}
      }
    \ifx#1\mathop\expandafter\displaylimits\fi}
\makeatother

\newcommand{\cupdot}{\charfusion[\mathbin]{\cup}{\cdot}}
\newcommand{\bigcupdot}{\charfusion[\mathop]{\bigcup}{\cdot}}

\usepackage[all]{xy}

\begin{document}

\title{ Homomorphisms of Fourier--Stieltjes Algebras  }
\author{Ross Stokke\footnote{
This research was partially supported   by an NSERC grant. }    }
\date{}
\maketitle

\begin{abstract}{\small   Every homomorphism $\vp: B(G) \ra B(H)$ between Fourier--Stieltjes algebras on locally compact groups $G$ and $H$ is determined by a continuous mapping $\alpha: Y \ra \DBG$, where $Y$ is a set in the open coset ring of $H$ and $\DBG$ is the Gelfand spectrum of $B(G)$ (a $*$-semigroup).  We exhibit a large collection of maps $\alpha$ for which $\vp=\ja: B(G) \ra B(H)$ is a completely positive/completely contractive/completely bounded homomorphism and establish converse statements in several instances. For example, we fully characterize all completely positive/completely contractive/completely bounded homomorphisms $\vp: B(G) \ra B(H)$ when $G$ is a Euclidean- or $p$-adic-motion group. In these cases, our description of the completely positive/completely contractive homomorphisms employs the notion of a  ``fusion map of a compatible system of homomorphisms/affine maps" and is quite different from the Fourier algebra situation. 

\smallskip

\noindent{\em Primary MSC codes:}  43A30, 43A22, 43A70,  47L25  

 }
\end{abstract}

Let $G$ and $H$ be locally compact groups.  The principal objects of study in abstract harmonic analysis are the  group and measure algebras $L^1(G)$ and $M(G)$ with convolution product and, dual to these algebras,  the Fourier and Fourier--Stieltjes algebras $A(G)$ and $B(G)$ with pointwise-defined product; $L^1(G)$ is a closed ideal in $M(G)$ and $A(G)$ is a closed ideal in $B(G)$.  When $G$ is abelian with dual group $\widehat{G}$,  $L^1(G)^\wedge$ is  $A(\widehat{G})$  and $M(G)^\wedge$ is  $B(\widehat{G})$, where $\mu \mapsto \hat{\mu}$ is the Fourier--Stieltjes transform on $M(G)$.  A fundamental endeavour in virtually every area of mathematics is to describe the structure-preserving maps between like objects.  The ``homomorphism problem" in abstract harmonic analysis asks for a description of all bounded algebra homomorphisms  of  $L^1(G)$ into   $M(H)$ and, dually,  of  $ A(G)$ into  $B(H)$.  After much effort, both versions of the problem were solved in the abelian case by Paul Cohen in 1960 \cite{Coh, Rud}.  For non-abelian groups,  both versions of the  problem---which are distinct in the non-abelian case---have been broadly studied by many authors, e.g., see \cite{Gre, Wal1, Il, Il-Sp, Il-Sp2,  Pham, Sto1, Sto}, but in full generality remain open.  In particular, Ilie and Spronk successfully generalized Cohen's results by describing all completely positive, completely contractive and completely bounded homomorphisms  $\vp: A(G) \ra B(H)$ when $G$ is amenable: any such $\vp$ is determined by a continuous map $\alpha:Y \ra G$ where $Y \in \Omega(H)$, the ring of sets generated by the open cosets of $H$, (we write $\vp = \ja$),  with $\alpha$ a homomorphism  and $Y$  a subgroup precisely when $\vp$ is completely positive; $\alpha$ an affine map  and $Y$  a coset precisely when $\vp$ is completely contractive; $\alpha$  a piecewise affine map  precisely when $\vp$ is completely bounded \cite[Theorem 3.7]{Il-Sp}.    

It seems equally natural to consider the problem of describing homomorphisms $\vp: B(G) \ra B(H)$. This problem has been settled when $\vp$ is an  isometric surjection \cite{Wal1},  $\vp$ is a positive/contractive isomorphism \cite{Pham},    $G$ is amenable and $\vp$ is a completely bounded isomorphism \cite{Pham},  or  $G$ is amenable and $\vp$ is a completely positive/completely contractive/completely bounded homomorphism that is   either weak$^*$- or strict-to-weak$^*$-continuous with respect to the strict topology on $B(G)$ determined by the ideal $A(G)$ \cite{Il-St}.  In each of these cases, $\vp$ takes the form---as with homomorphisms $\vp: A(G) \ra B(H)$---$\vp=\ja$ where $\alpha: Y \ra G$ is a continuous mapping on a set $Y$ in $\Omega(H)$;  in the language introduced in Section 1, any such $\vp$ is ``spatial".  However, not every homomorphism $\vp: B(G) \ra B(H)$ is spatial, though very little is known about non-spatial homomorphisms; the only ad hoc constructions of such homomorphisms seem to be found in two proofs in  \cite{Rud} and \cite{Il-St}.  In fact, beyond the results from \cite{Wal1,  Il-St, Pham}, very little is known about the homomorphisms  $\vp: B(G) \ra B(H)$.       
In light of Proposition \ref{FSAlgHomGeneralProp}, which shows that every homomorphism $\vp:B(G) \ra B(H)$ is determined by a continuous map  $\alpha:Y \ra \DBG$ where  $Y$ in $\Omega(H)$ and $\DBG$ is the Gelfand spectrum of $B(G)$,  the  gap in our understanding of Fourier--Stieltjes algebra homomorphisms may be attributed, in part, to the inaccessibility of $\DBG$ in comparison with $\Delta(A(G))$, which is just $G$ itself.  To quote  E. Kaniuth and A.T.-M. Lau, ``when $G$ is a noncompact locally compact abelian group, according to common understanding, the spectrum of $B(G)= M(\widehat{G})$ is an intractable object" \cite[Section 2.9]{Kan-Lau}.  Nevertheless, the structure of $\DBG$ and $\Delta(A)$  for certain subalgebras $A$ of $B(G)$ has been successfully studied by various authors, e.g., \cite{Dun-Ram, Wal2, Liu-Mis, Il-Sp2},  and shown to be fairly accessible in some---typically non-abelian---cases.   Results from \cite{Wal2} and \cite{Il-Sp2}, and the work herein, allow us to address our  main questions, which  are as follows:    

 \bi  \item[(i)] Given a set $Y \in \Omega(H)$,  can we systematically construct  a large collection of (non-spatial) continuous maps $\alpha: Y \ra  \DBG$  that determine homomorphisms $\vp = j_\alpha:B(G) \ra B(H)$, and  what properties are possessed by these homomorphisms?  
 
 \item[(ii)]  Given a class of homomorphisms---e.g., homomorphisms $\vp = \ja: B(G) \ra B(H)$  that are (completely) positive/ (completely) contractive/ completely bounded---does $\alpha: Y\ra \DBG$ take the form  one might expect based on  the  existing theory of homomorphisms from $A(G)$ into $B(H)$?  
 \ei 
 
We show that $\vp= \ja:B(G) \ra B(H) $ is a completely positive homomorphism when $\alpha$ is a continuous homomorphism of an open subgroup of $H$ into a  subgroup of $\DBG$ (Corollary \ref{CPHomCor}); a completely contractive homomorphism when $\alpha$ is a continuous affine map into a  subgroup of $\DBG$ (Proposition \ref{CCProp}); and  a completely bounded homomorphism when $\alpha$ is what we shall call  a continuous  pw$^2$-affine map (Proposition \ref{PWAffineCBProp}). In this way,  we provide a positive answer to question (i) and  show that spatial homomorphisms are the exception rather than the norm.  
We show that the answer to question (ii) is sometimes, but not always.  When $\alpha$ maps into a subgroup of $\DBG$, we show that $\vp=\ja$ is completely positive/completely contractive if and only if $\alpha$ is a homomorphism/affine map (Corollaries  \ref{CPCor1} and \ref{CC Thm Cor}).  In light of the theorems of Cohen and Ilie--Spronk, it is perhaps  tempting to conjecture that $\alpha$ must always be a  homomorphism/affine map into a subgroup of $\DBG$ when $\vp = \ja$ is completely positive/completely contractive. However, we show that  when $\alpha$ maps into the union of a certain type of  compatible system of subgroups of $\DBG$, $\vp = \ja$ is completely positive/completely contractive if and only if $\alpha$ is what we shall call a ``fusion map of a compatible system of homomorphisms/affine maps" (Theorems \ref{Second Main Thm} and \ref{CCThm}) and show that it is easy to construct examples of these fusion maps (Example \ref{Construction of compatible homoms Ex}). (Note that we do not require  $G$ to be amenable in any of these general results.)   In particular, we  characterize all completely bounded, completely contractive and completely positive homomorphisms $\vp: B(G) \ra B(H)$ when $G$ is  a Euclidean- or $p$-adic motion group (Theorem \ref{AFG Main Homom Thm} and  Corollary \ref{Full characterization of cc and cp in special cases Cor}). More generally, we will consider homomorphisms $\vp: A \ra B(H)$ where $A$ is any closed translation-invariant unital subalgebra of $B(G)$, with applications to some specific examples of the spine of $B(G)$  (studied by Ilie and Spronk in \cite{Il-Sp2}).  Along the way, we extend several results concerning $\DBG$ from \cite{Wal2}  with simple new proofs.   The basic theory of multiplicative domains of completely positive operators, the origin of which is \cite{Choi} by M.-D. Choi, is a fundamental tool in Section 4.

We proceed with a brief review of relevant terms and results.   Throughout this paper, $G$ and $H$ are locally compact groups.  We refer the reader to  \cite{Ars, Eym, Kan-Lau} for definitions and basic properties of the Fourier algebra $A(G)$ and the Fourier--Stieltjes algebra $B(G)$; unless stated otherwise, we follow the notation used in these references.     By a  representation $\{ \pi, \H\}$  of $G$ we  will always mean a continuous unitary representation, $\pi$,  on a Hilbert space $\H$; ${\cal L}(\H)$ denotes the space of  bounded linear operators on $\H$.   The closed linear span in $B(G)$  of all coefficient functions   $$ \xi *_\pi \eta(s) := \l \pi (s) \xi | \eta \r \qquad  (s \in
G, \ \xi, \eta \in \H)$$ is denoted $A_\pi$. The von Neumann subalgebra of ${\cal L}(\Hp)$ generated by $\pi(G)$, denoted $VN_\pi$, is identified with the dual of $A_\pi$ through the pairing 
$\l T, \xpe \r = \l T \xi | \eta\r_\H$;  hence,  $\l \pi(s), u \r = u(s)$ for $s \in G$ and $u \in A_\pi$.  The Fourier spaces $A_\pi$ are precisely the closed translation-invariant subspaces of $B(G)$. Letting $P(G)$ denote the set of continuous positive definite functions on $G$, $A_\pi \cap P(G)$ is the set of normal positive linear functionals on $VN_\pi$. Letting $\{ \omega_G, \H_{\omega_G}\}$   and $\{ \lambda_G, L^2(G)\}$  respectively denote the universal and left regular representations of $G$, $B(G) = A_{\omega_G}$ and $A(G) = A_{\lambda_G}$; as usual, we write $W^*(G)$ and $VN(G)$  for $VN_{\omega_G}$  and $VN_{\lambda_G}$, respectively.    

As the predual of a von Neumann algebra, $A_\pi$ has a canonical operator space structure that  agrees with the subspace operator space structure inherited from $B(G) = W^*(G)_*$.   Moreover, $A_\pi$ is a completely contractive Banach algebra when it is also a  subalgebra of $B(G)$.  Over the course of the last 25 years, studying the Fourier and Fourier--Stieltjes algebras as completely contractive Banach algebras has proven invaluable to the theory, e.g., see \cite{Run2, Spr}.  We say that a map $\vp: A _\pi \ra B(H)$ is (completely) positive when its dual map $\vp^*: W^*(H) \ra VN_\pi$ is (completely) positive.   Equivalently, $\vp: A_\pi  \ra B(H)$ is positive if and only if $\vp$ maps positive definite functions in $A_\pi$ to positive definite functions in $B(H)$. We follow the standard references  \cite{Eff-Rua} and   \cite{Pau} on the theory of operator spaces and completely bounded maps.

   By a topological isomorphism of  topological groups, we shall mean a group isomorphism that is also a homeomorphism.  As noted in \cite{Il}, a subset $E$ of $H$ is a coset of some subgroup of $H$ exactly when $E E^{-1}E = E$, and a map $\alpha: E \rightarrow G$ is called \it affine \rm if for any
$x,y,z \in E$, $\alpha(xy^{-1}z) =
\alpha(x)\alpha(y)^{-1}\alpha(z)$. Note that  for any $y_0 \in E$, $H_0 = y_0^{-1} E=E^{-1}E$ is a subgroup of $H$,  and the map defined by $\beta(h) = \alpha(y_0)^{-1}
\alpha(y_0h)$ $(h \in H_0)$ is a homomorphism of $H_0$ into $G$ when   $\alpha$ is an
affine map; conversely, if $\beta:H_0 \ra G$ is a homomorphism and $x_0\in G$, then $\alpha(x)= x_0\beta(y_0^{-1}x)$ defines an affine map on $E$ (see
\cite[Remark 2.2]{Il}). Thus, affine maps are exactly the translates of subgroup homomorphisms.  Letting $\Omega(H)$ denote the ring of sets generated by the open cosets of $H$,  every set in $\Omega(H)$ can be expressed as a finite union of disjoint sets in
$$
\Omega_0(H)=\left\{E_0 \setminus \left(\bigcup_{1}^{m}E_k\right):\begin{array}
{l}E_0 \subseteq H
{\rm \ an \  open \ coset,}\\
              E_1,\ldots,E_m {\rm\  open \ subcosets\ of \ infinite \ index \ in}\
E_0 \end{array}\right \}.
$$  A map $\alpha : Y \rightarrow G$ is \it
piecewise affine \rm if
 there exist pairwise disjoint sets $Y_1,...,Y_n \in \Omega_0(H)$ such that $Y=\bigcup_{i=1}^{n}{Y_i}$ and for each $i$,  $\alpha\big{|}_{Y_i}$ has an affine
  extension $\alpha_i$ mapping $Aff(Y_i)$, the smallest coset containing $Y_i$,  into $G$ \cite{Coh, Il}. \rm

We let $\Delta(A)$  denote the Gelfand spectrum of a commutative Banach algebra $A$.


  \section{Preliminary results} 
  
  \subsection{Completely positive maps} 
  
  We will make repeated use of a few basic facts about completely positive maps, collected below for the convenience of the reader.  References are, for example, 
 \cite[Ch. 2, 3 and 4]{Pau} and \cite[Ch. 5]{Eff-Rua}.   Regarding part (c), note that either hypothesis implies that $\psi_2$ satisfies the Schwarz inequality \cite[Proposition 3.3]{Pau}, which is the essential ingredient in the proofs of  \cite[Theorem 3.18]{Pau} and \cite[Corollary 5.2.2]{Eff-Rua}. 

  \bt  \label{CPBasicThm} Let $A$ and $B$ be  $C^*$-algebras with identities, $\psi: A \ra B$ a  linear map. 
  \bi \item[(a)] If $\psi$ is (completely) positive, then $(\|\psi\|_{cb}=) \| \psi \| = \| \psi(1_A) \|$ and for each $a \in A$, $\psi(a^*) = \psi(a)^*$.
  \item[(b)]  If $\psi$ is (completely)  contractive and unit preserving, then $\psi$ is (completely) positive. 
  \item[(c)]  Suppose that $\psi$ is either completely positive and (completely) contractive, or 4-positive and unit preserving. If $\psi(b)^* \psi(b) = \psi(b^*b)$ $(\psi(b) \psi(b)^* = \psi(bb^*))$ for some $b \in A$, then $\psi(ab) = \psi(a)\psi(b)$ $(\psi(ba) = \psi(b) \psi(a))$ for every $a \in A$.  
  \ei   \et


  As an application of Theorem \ref{CPBasicThm}, we now offer a very  short new proof of a theorem due to Ilie--Spronk \cite{Il-Sp} and Pham \cite{Pham}.  The proof of Theorem \ref{CPMainThm} contains similar elements, but this proof requires no additional machinery and is likely of independent interest.

  Let $\vp: A(G) \ra B(H)$ be a homomorphism, $\vp^*: W^*(H) \ra VN(G)$ the dual map.  Then for each $h \in H$, $\vp^*(h) \in \Delta(A(G)) \cup \{0\} = G \cup \{0\}$,  (identifying $G$ with $\lambda_G(G)$). Letting $Y= H \bs \ker \vp^*$, $Y$ is an open subset of $H$ and $\alpha: Y \ra G: h \ra \vp^*(h) $ is a continuous map such that $$\displaystyle{\vp(u)(h) = \l \vp^*(h), u \r = \left \{ \begin{array}{ll}
                          u(\alpha(h))  & \mbox{ $ h \in  Y$}\\
                          \ \ \  0 & \mbox{  $ h \in H \bs Y$}
                            \end{array}
                 \right.  = \ja u(h).}   $$

 \bp \label{PhamProp5.8} Suppose that $\vp: A(G) \ra B(H)$ is a nontrivial 4-positive homomorphism. Then there is an open subgroup $H_0$ of $H$ and a continuous homomorphism  $\alpha: H_0 \ra G$  such that $\vp = \ja$.  
 \ep 
 
 \begin{proof}   Letting $\psi = \vp^*$, $0 \neq \|\vp\| = \|\psi\| = \|\psi(e_H)\|$, so $e_H \in Y$. Let $s, t \in Y$. Then $\psi(t^{-1}) = \psi(t^*) = \psi(t)^* = \psi(t)^{-1} \neq 0$, so $t^{-1}\in Y$. Since $\psi(e_H) = \psi(e_H)^{-1}$,  $\psi(e_H)^2 = e_G$ and $\psi(e_H) \geq 0$; hence $\psi(e_H) = (e_G)^{1/2} = e_G$. It follows that 
 $\psi(t)^* \psi(t) = \psi(t)^{-1} \psi(t) = e_G = \psi(t^* t)$ and therefore $\psi(st) = \psi(s) \psi(t) \neq 0$. Hence, $st \in Y$ and $\alpha(st) = \alpha(s) \alpha(t)$.  \end{proof}

 \br   \label{PhamProp 5.8 Remark}  \rm 1. Ilie and Spronk proved Proposition \ref{PhamProp5.8} as one part of \cite[Theorem 3.7]{Il-Sp} under the slightly stronger assumption that $G$ is amenable and $\vp$ is completely positive, while Pham proved it under the slightly weaker assumption that $\vp$ is 2-positive \cite[Proposition 5.8]{Pham}. \\
 2. By translation, one can now  use Theorem \ref{CPBasicThm}(b) and the argument found in the short proof of \cite[Theorem 5.1]{Pham} to show that if $\vp: A(G) \ra B(H)$ is a 4-contractive homomorphism, then  $\vp = \ja$ for some continuous affine map $\alpha:Y \ra G$ where $Y$ is an open coset of $H$. \\  
 3. For future reference, observe that the arguments found in the proof of  Proposition \ref{PhamProp5.8} and the above remark  also provide characterizations of all 4-positive and 4-contractive homomorphisms $\vp: A \ra B(H)$ whenever $A$ is a closed translation-invariant subalgebra of $B(G)$ for which $\Delta(A) = G$.  For example,  when $A = B_0(G)$, the Rajchman algebra of $G$,  is regular, $\Delta(A) = G$ \cite[Corollary 2.2]{Kan-Lau-Ulg}; this is the case, when $G$ is a connected semisimple Lie group with finite centre \cite[Example 2.6(3)]{Kan-Lau-Ulg} or $G$ is the Euclidean motion group  $\mathbb{R}^n \rtimes SO(n)$ and $n \geq2$  \cite[Theorem 4.1]{Kan-Lau-Ulg}.
 \er

    We can now give a new proof of \cite[Theorem 2.1]{Il-Sp} due to Ilie and Spronk.  This should be of independent interest since it provides another example of the use of operator space theory in abstract harmonic analysis \cite{For-Woo, Run2, Spr} and makes this article largely self-contained.   
  
  \bc \label{Il-SpThm2.1Cor}  Let $u$ be a non-zero idempotent in $B(G)$. 
  \bi \item[(a)] The following statements are equivalent: 
  \bi  \item[(i)] $\|u \|=1$ and $u(e_G) \neq0$;
  \item[(ii)]  $u = 1_H$ for some open subgroup $H$ of $G$;
  \item[(iii)] $u$ is positive definite.
      \ei 
  \item[(b)] $\|u \|=1$ if and only if $u=1_C$ for some open coset $C$ in $G$.
  \ei 
  \ec  
  
  \begin{proof} As $u^2=u$ and $u$ is continuous, $u = 1_Y$ for some open and closed subset $Y$ of $G$. As noted -- and proved -- in \cite{Il-Sp}, (ii) implies (iii) and sufficiency in (b) are well-known facts. 
  
  (a) Assume condition (i) holds.   Since $1_Y(e_G)  \neq 0$, $e_G \in Y$ and, because $B(G)$ is a completely contractive Banach algebra and $\|1_Y\|=1$,   $\vp(v)=v1_Y$ defines a completely contractive homomorphism $\vp: A(G) \ra B(G)$. The dual map $ \vp^*: W^*(G) \ra VN(G)$ is hence a complete contraction satisfying $\vp^*(s) =s$ if $s \in Y$ and $\vp^*(s) = 0 $ if $s \in G \bs Y$. 
  Hence, $\vp = \ja$ where $\alpha: Y \ra G: s \mapsto s$ and, as a unit-preserving complete contraction, $\vp^*$  is  completely positive.   That is, $\ja$ is completely positive so $Y$ is an open subgroup by Proposition \ref{PhamProp5.8}.  Hence (i) implies (ii), and we have already noted that (ii) implies (iii).  For (iii) implies (i), observe that $ \|u\| = u(e_G) = 1_Y(e_G)$ and $u$ is nonzero, so $u(e_G) \neq 0$ and  $\|u\|= 1$. 
  
  (b)  Suppose that $\|u \| = \|1_Y\| =1$. Taking $y \in Y$, $\ell_y 1_Y = 1_{y^{-1}Y} $ is also a norm-one idempotent in $B(G)$, (e.g., by \cite[(2.19)]{Eym} or, more generally, Lemma \ref{TranslationCCHomLem}).  Since $\ell_y 1_Y(e_G) = 1$, $y^{-1}Y $ is an open subgroup of $G$ by part (a).    
  \end{proof}

  \subsection{The spectrum  of $B(G)$ and its closed subalgebras}
  
  The following Arens-product formulation of the product in $VN_\pi= A_\pi^*$ will be often be  used  in this paper. It should be known---e.g., it is implicit in \cite{Spr-Sto}---but we provide the simple proof below for the convenience of the reader.  Observe that if $\tilde{u}(s) = \ov{u(s^{-1})}$, then $(\xpe)^\sim = \eta *_\pi \xi$; hence, $\tilde{u} \in A_\pi$ whenever $u \in A_\pi$.  
    
  \bp  \label{ArensProductProdProp}  Let $m,n \in VN_\pi$, $u \in A_\pi$, $s,t \in G$. Letting 
  $$u\c s(t) = u(st), \ \ s \c u(t) = u(ts), \ \  n\c u(s) = n(u \c s) \ \ {\rm and} \ \ u \c m(s) = m(s \c u), $$
the product and involution in $VN_\pi$ are given by   
$$\l mn, u \r= \l m, n\c u\r = \l n , u \c m \r \quad {\rm and} \quad \l m^*, u \r = \overline{\l m, \tilde{u}\r}.$$
  \ep 

\begin{proof} Since $VN_\pi$ is a dual Banach algebra, $A_\pi$ is a $VN_\pi$-submodule of $VN_\pi^*$ via 
$$ \l m, n\c u\r = \l mn, u \r = \l n, u \c m \r$$
\cite{Run}.  But $$\pi(s) \c u (t) = \l \pi(t), \pi(s) \c u \r = \l \pi(ts), u \r = u(ts) = s \c u(t),$$so 
$$u \c m(s) = \l \pi(s), u \c m \r = \l m, \pi(s) \c u \r = \l m ,  s \c u \r;$$ similarly,  $n \c u(s) = \l n, u\c s\r$. If $m = \sum \lambda_i \pi(s_i)  \in {\rm span}\{ \pi(G)\}$, then 
$$m^*(u) = \l \sum \ov{\lambda_i} \pi(s_i^{-1}), u \r = \sum  \ov{\lambda_i} u(s_i^{-1})  = \ov{m(\tilde{u})}.$$ Weak$^*$-density of  ${\rm span}\{ \pi(G)\}$ in $VN_\pi$ and weak$^*$-continuity of involution now yields the formula for any $m \in VN_\pi$. \end{proof} 

 Walter gave a different proof of the following result when  $A= B(G) = A_{\omega_G}$ in \cite{Wal2}. 
 
\bc    \label{SpectrumSemigpCor} Let $A$ be a  closed unital translation-invariant subalgebra of $B(G)$.  Then the spectrum of $A$, $\Delta(A)$, is a compact semitopological $*$-semigroup. \ec
 
 \begin{proof}  We know that $A= A_\pi$ for some $\pi$, so $\Delta(A)$ is contained in $VN_\pi$.    Since $1_G \in A$, $\Delta(A)$ is compact.  Let $m,n \in \Delta(A)$, $u,v \in A$. Then for $s \in G$, 
 $$n \c (uv)(s) = \l n , (uv) \c s \r = \l n, (u \c s)(v \c s) \r = \l n, u\c s\r \l n, v \c s \r = (n \c u)(s) (n \c v)(s)$$ so 
 $$ \l mn, uv \r = \l m, n \c (uv) \r = \l m, (n \c u)(n \c v) \r =  \l m, n \c u \r \l m, n \c v \r  = \l mn, u\r \l mn, v\r.$$
Hence, $mn \in  \Delta(A) \cup \{ 0 \}$. However, $(n \c 1_G)(s) = n (1_G \c s) = n(1_G) = 1$, so $mn(1_G) = m(n \c 1_G) = m(1_G) =1$; hence $mn \in \Delta(A)$.   Using the formula for $m^*$ from Proposition \ref{ArensProductProdProp}, it is also easy to show that $m^* \in \Delta(A)$.  Since multiplication is separately continuous in $VN_\pi$, $\Delta(A)$ is semitopological.  \end{proof} 

 
  \bi \item[] \it Throughout the remainder of this subsection, we assume that  $A=A_\pi$ is a closed translation-invariant subalgebra of $B(G)$  for which $\Delta(A)$ is a $*$-semigroup. \rm \ei

For example,  this is the case  when $A$ is    unital    (Corollary \ref{SpectrumSemigpCor}) and  when  $A=A(G)$.  Suppose further that $A=A_\pi$ is point separating. Then $\pi$  is a continuous group isomorphism of $G$ onto its image $\pi(G)$ in $\Delta(A) \subseteq (VN_\pi)_{\| \c \| =1}$, and we shall often identify $G$ and $\pi(G)$ as groups. If $A$ is a regular algebra of functions on $G$---e.g., if  $A$ contains $A(G)$---then  it is easy to see that  $\pi$ is actually a topological isomorphism, and we can identify $G$ and $\pi(G)$ as topological groups. 

  The Arens-product approach also allows us to give an easy proof of another useful result;  in the case that $A=B(G)$, a different proof is   given by Walter  in  \cite{Wal1}.   
 
  \bc  \label{Wal1Cor}  Suppose that $A=A_\pi$ contains $A(G)$. The following statements hold: 
 \bi \item[(i)] $A(G)^\perp$ is a weak$^*$-closed ideal in $VN_\pi$.
 \item[(ii)] $\Delta(A) \backslash G = \Delta(A) \cap A(G)^\perp = h(A(G))$, the hull of $A(G)$; hence $\Delta(A) \backslash G$ is a closed ideal in $\Delta(A)$.
 \item[(iii)]  $\Delta(A) \cap (VN_\pi)_r = \Delta(A)\cap (VN_\pi)_u = G$, where $(VN_\pi)_r$ and $(VN_\pi)_u$ respectively denote the invertible and unitary elements in $VN_\pi$. 
 \ei 
 \ec 
 
 \begin{proof} (i) More generally, letting $B$ be any translation-invariant subset of $A$, we observe that $B^\perp$ is a weak$^*$-closed ideal in $VN_\pi$:   Let $m \in VN_\pi$, $n \in B^\perp$. For $v \in B$ and $s \in G$, $n\c v(s) = n ( v\c s) = 0$ and $v \c n(s) = n(s \c v) = 0$, so  $mn(v) = m (n \c v) = 0$ and $nm(v) = m (v \c n) =0$.  (This is a special case of a standard result, e.g.,  \cite[Theorem III.2.7]{Tak}.)
  
 (ii)  From (i), $h(B) = \Delta(A) \cap B^\perp$ is a closed ideal in $\Delta(A)$ whenever $B$ is a translation-invariant subset of $A$, so  we only need to  establish the first set equality in (ii). Since $G$ is disjoint from $A(G)^\perp$,  the first set contains the second one. To establish the reverse containment, supposing  that $m \in \Delta(A) \bs A(G)^\perp$,  we show $m \in G$ via the argument used to establish the lemma on p. 27 of \cite{Wal1}: $m\big{|}_{A(G)} \in \Delta(A(G))=G$, so $m\big{|}_{A(G)} = g$ for some $g \in G$. Taking $v \in A(G)$ such that $m(v) = v(g) =1$, for each $u \in A$, $$m(u) = m(u) m(v) = m(uv) = uv(g) = u(g),$$ since $A(G)$ is an ideal in $A$. Hence $m = \pi(g) \in G$. 
 
 (iii) Given $m \in   \Delta(A) \cap (VN_\pi)_r$,  take $n \in VN_\pi$ such that $mn = e_G$. If  $m \notin G$, then $m \in A(G)^\perp$ by (ii), whence $e_G = mn \in A(G)^\perp$ by (i), (a contradiction). Hence $m \in G$. Since elements of $G(=\pi(G))$ are unitary in $VN_\pi$ and unitary operators are invertible, we are done. 
 \end{proof} 
 
 \subsection{Homomorphisms of Fourier-Stieltjes algebras: introduction} 
 
 Let $\Omega(H)$ denote the ring of sets generated by the cosets of open subgroups of $H$.  
 
  \bp \label{FSAlgHomGeneralProp} Let $A= A_\pi$ be a closed translation-invariant subalgebra of $B(G)$.  For a linear map $\vp: A \ra B(H)$, the following statements are equivalent:  
 \bi \item[(i)]  $\vp$ is a homomorphism; 
 \item[(ii)] $\vp $ is bounded and $\vp^*$ maps $H$ into $\Delta(A) \cup \{0\}$ where $\vp^*: W^*(H) \ra VN_\pi$ is the dual map of $\vp$;
 \item[(iii)]  there is an open  subset $Y$ of $H$ and a continuous map $\alpha:Y \ra \Delta(A)$ such that 
$$\displaystyle{\vp(u)(h) = \ja(u)(h)\coloneqq \left \{ \begin{array}{ll}
                           \l \alpha(h), u \r   & \mbox{ $ h \in  Y$}\\
                          \ \ \ \ 0 & \mbox{  $ h \in H \bs Y.$}
                            \end{array}
                 \right. } $$
 \ei 
 When $\vp$ is a homomorphism,  $Y$ and $\alpha$ in (iii)  are uniquely determined by $\vp$  as follows: $Y= H \bs \ker (\vp^*)$ and $\alpha = \vp^* \large{|}_Y$. If $A$ is unital, then $Y \in \Omega(H)$. 
 \ep 
 
 \begin{proof}  Assume (i) holds. Since $B(H)$ is semisimple, $\vp$ is bounded \cite[Corollary 2.1.10]{Kan},  and for each $h \in H$, 
 $$\vp^*(h)(uv) = \vp(uv)(h) = \vp(u)(h)\vp(v)(h) = \vp^*(h)(u)\vp^*(h)(v),$$
so $\vp^*(H)$ is contained in  $\Delta(A) \cup \{0\}$.    Suppose that (ii) holds. Letting $Y= H \bs \ker (\vp^*) =  \{ h \in H: \vp^*(h) \neq 0\}$, $Y$ is open, $\alpha: Y \ra \Delta(A): y \ra \vp^*(y)$ is a continuous map and $\vp = \ja$ as defined in (iii). Suppose $1_G \in A$. Since $\l \gamma, 1_G\r = 1$ for $\gamma \in \Delta(A)$, $\vp(1_G) = \ja(1_G) = 1_Y$, an idempotent in $B(H)$. By Host's idempotent theorem \cite{Hos}, $Y \in \Omega(H)$. Finally, assume that (iii) holds and define $\alpha_0: H \ra \Delta(A) \cup \{0\}$ by putting $\alpha_0\big{|}_Y = \alpha$ and  $\alpha_0\big{|}_{H \bs Y} = 0$. For $h \in H$,  $\vp(u)(h) = \ja u(h)  = \l \alpha_0(h), u \r$, so 
$$\vp(uv)(h) = \l \alpha_0(h), uv\r = \l \alpha_0(h), u \r \l \alpha_0(h), v\r = \vp(u)(h) \vp (v)(h). $$
Hence,  $\vp$ is a homomorphism.  It is easy to see that $Y$ and $\alpha$ are uniquely determined by $\vp$.  \end{proof} 
 
 
 
 

\noindent {\bf Notation:}    \it  Unless stated otherwise, whenever $\vp:A \ra B(H)$ is a homomorphism, we will use the notation $\vp = \ja$ where $\alpha:Y \ra \Delta(A)$ is a continuous mapping on an open subset $Y$ of $H$ and   $\ja$ is defined as in Proposition \ref{FSAlgHomGeneralProp}(iii). \rm

\medskip

 It is helpful to introduce some new terminology. We will say that the homomorphisms $\vp, \psi: A \ra B(H)$ are \it disjoint \rm if  $\vp(u) \psi(w) = 0$ for every $u, w \in A$. If $C$ and $D$ are disjoint sets and $f:C \ra X$, $g:D\ra X$, then $f\cupdot g: C\cupdot D \ra X$ is the map that agrees with $f$ on $C$ and $g$ on $D$.
When $A $ is point separating, we will say that $\vp= \ja$ is \it spatial \rm if $\alpha(Y) \subseteq G$ and \it singular \rm if $\alpha(Y) \subseteq \Delta(A) \bs G$.   The following useful lemma clarifies these notions. 

\blem \label{Disjoint Homom Lemma}   Let $\vp = j_\alpha, \psi= j_\beta : A \ra B(H)$ be homomorphisms, where  $\alpha:Y \ra \Delta(A)$ and  $\beta : Z \ra \Delta(A)$. Then $\vp $ and $\psi$ are disjoint if and only if $Y \cap Z = \emptyset$. Moreover, $\vp + \psi:A \ra B(H)$ is a homomorphism whenever $\vp$ and $\psi$ are disjoint and $\vp + \psi = j_{\alpha\cupdot\beta}$.  
\elem 
 
 \begin{proof}  This is readily checked. \end{proof} 
 
 \bi \item[] \it Throughout the remainder of this section, we assume that  $A=A_\pi$  is a closed unital translation-invariant subalgebra of $B(G)$ that contains $A(G)$.  \rm \ei 
 
 We  assume  $A$ is unital to ensure that $\Delta(A)$ is a $*$-semigroup and $Y \in \Omega(H)$.  However, most results in this section  hold whenever $\Delta(A)$ is known to be a $*$-semigroup and  $Y \in \Omega(H)$.  
 
  Our assumption implies that there is a continuous unitary representation $\sigma$ of $G$ that is disjoint from $\lambda_G$ such that $\pi = \lambda_G \oplus \sigma$ and   $A = A(G) \oplus_1 A^s(G)$, where $A^s(G) = A_\sigma$ is what we shall call the  \it singular  part \rm  of $A$   \cite[Proposition 3.12, Corollary 3.13 and Theorem 3.18]{Ars}.    When $A=  B(G)$, this is the \it Lebesgue decomposition \rm of $B(G)$ that was studied in \cite{Miao} and \cite{Kan-Lau-Sch}.   With the following lemmas, we record more useful facts for future reference.  

\blem  \label{Basic Lemma spatial singular disjoint}  Let $\vp = j_\alpha: A \ra B(H)$ be a homomorphism. 
\bi \item[(i)] $\vp$ is spatial if and only if for each $y \in Y$, $\vp(v)(y) \neq 0$ for some $v \in A(G)$. 
\item[(ii)]  The following statements are equivalent: 
\bi \item[(a)]  $\vp$ is singular;
\item[(b)] $A(G) $ is contained in the kernel of $\vp$;
\item[(c)] for each $w = v + u \in A = A(G) \oplus_1 A^s(G)$, $\vp (w) = \vp (u)$.
\ei 
\ei 
\elem 
\begin{proof} Using Corollary \ref{Wal1Cor}(ii), this is easily checked.  \end{proof}
 
 \blem \label{Restriction to A(G) Lemma}   Let $\vp = j_\alpha: A \ra B(H)$ be a homomorphism, $\vp_F:A(G) \ra B(H)$ the restriction of $\vp$ to the Fourier algebra $A(G)$. Then $Y_a = \alpha^{-1}(G) \subseteq Y$ is an open subset of $H$ such that $\vp_F = j_{\alpha_a}$ where $\alpha_a: Y_a \ra G$ is the restriction of $\alpha$ to $Y_a$. 
 \elem
 
 \begin{proof} This is an easy consequence of  Corollary \ref{Wal1Cor}(ii).  
 \end{proof} 
 
 \bp \label{Lebesgue decomposition Prop} Let  $\vp = j_\alpha: A \ra B(H)$ be a (completely bounded)  homomorphism, $Y_a = \alpha^{-1}(G)$,  $\alpha_a: Y_a \ra G$ the restriction of $\alpha$ to $Y_a$.  \bi  \item[(i)] The following statements are equivalent: \bi \item[(a)]  there are disjoint  spatial and singular (completely bounded)  homomorphisms, $\vp_a: A \ra B(H)$ and $\vp_s:A \ra B(H)$ respectively, such that $\vp = \vp_a + \vp_s$; 
 \item[(b)] $Y_a \in \Omega(H)$;  
 \item[(c)]  $j_{\alpha_a}$ defines a (completely bounded) spatial homomorphism from $A$ into $B(H)$. 
 \ei 
   \item[(ii)] When the decomposition from statement (i)(a)  exists,  $\vp_a$ and $\vp_s$ are unique,  $\vp_a = j_{\alpha_a}$ where $\alpha_a :Y_a \ra G$ is a continuous map on $Y_a \in \Omega(H)$,  $\vp_s = j_{\alpha_s}$ where $\alpha_s:Y_s \ra \Delta(A)\bs G$ is a continuous map on $Y_s \in \Omega(H)$,  and $\alpha_a \cupdot \alpha_s = \alpha$.  
 \ei 
 \ep 
 
 \begin{proof}  (i) Assuming (c), (b) follows by Proposition \ref{FSAlgHomGeneralProp} since $A$ is unital.   Assume (b). Let $Y_s = Y\bs Y_a \in \Omega(H)$, $\alpha_s = \alpha{\big |}_{Y_s}$.   Then $1_{Y_a}$ and $1_{Y_s}$ are idempotents in $B(H)$, so $\vp_a, \vp_s: A \ra B(H)$ defined by $\vp_a(u) = \vp(u) 1_{Y_a}$, $\vp_s(u) = \vp(u) 1_{Y_s}$  are (completely bounded) homomorphisms and it is easy to check that $\vp_a = j_{\alpha_a}$ and $\vp_s = j_{\alpha_s}$.  Hence, $\vp_s$ and $\vp_a$ are disjoint spatial and singular homomorphisms such that $\vp_a + \vp_s = j_{\alpha_a \cupdot \alpha_s} = j_\alpha = \vp$. This is (a).  
 Suppose (a) holds. Let $\vp_a = j_\gamma$ and $\vp_s = j_\beta$ where $\gamma:Y_\gamma \ra G$ and $\beta: Y_\beta \ra \Delta(A)\bs G$. By Lemma \ref{Disjoint Homom Lemma}, $j_\alpha = \vp = j_{\gamma \cupdot \beta}$, so---by the uniqueness part of  Proposition \ref{FSAlgHomGeneralProp} ---$\alpha= \gamma \cupdot \beta$. Hence, $Y_a = \alpha^{-1}(G) = Y_\gamma$ and $\alpha_a = \alpha{\big |}_{Y_a} = \gamma$. Thus, $j_{\alpha_a} = j_\gamma = \vp_a$, a  (completely bounded) homomorphism of $A$ into $B(H)$; hence, (c) holds.  
  
  Only  the statement about uniqueness in (ii) still requires justification. Suppose we have another decomposition $\vp = \vp_1 + \vp_2$ where $\vp_1 = j_{\gamma}$, $\gamma: Y_\gamma \ra G$, $\vp_2 = j_{\beta}$, $\beta: Y_\beta \ra  \Delta(A)\bs G$, are disjoint spatial and singular homomorphisms respectively.  Then $Y_\gamma$ and $Y_\beta$ are disjoint and $\vp = j_{\alpha_a\cupdot \alpha_s} = j_{\gamma \cupdot \beta}$, so the uniqueness part of Proposition \ref{FSAlgHomGeneralProp} yields \beqs  \label{Spatial Singular Decomp Eqn}  Y=Y_a \cupdot Y_s = Y_\gamma \cupdot Y_\beta \quad {\rm and } \quad  \alpha = \alpha_a \cupdot 
 \alpha_s = \gamma \cupdot \beta.   \eeqs
 The argument establishing  (a) implies (c) in (i) shows that $ Y_a = Y_\gamma$ and $ \alpha_a = \gamma$.   
 We conclude  that  $Y_s= Y_\beta$ and $\alpha_s = \beta$.   \end{proof} 
  
  For the sake of brevity, when $\vp: A \ra B(H)$ has a decomposition $\vp = \vp_a + \vp_s$ as in Proposition \ref{Lebesgue decomposition Prop}, we shall refer to it as its \it Lebesgue decomposition. \rm When it exists, $\vp_a$  is the unique  spatial extension   of  $\vp_F$, as defined in Lemma \ref{Restriction to A(G) Lemma}, to $A$.

 When $\vp_a = j_{\alpha_a}:A \ra B(H)$ is spatial, its restriction to $A(G)$ is $j_{\alpha_a} :A(G) \ra B(H)$, so  $\alpha_a$ has been described in many instances \cite{Coh, Il-Sp, Pham}. Thus, when $\vp$ has a Lebesgue decomposition $\vp = \vp_a + \vp_s$,  we can study its spatial and singular parts separately and can often provide a detailed description of  its spatial part $\vp_a$ (Corollary \ref{Lebesgue Decomp Cor}).  As we shall see,  we can sometimes provide a detailed description of its singular part $\vp_s$ as well. We do not know if $\vp$ always has a Lebesgue decomposition, but will provide an affirmative answer in several cases below.   
  
  
  
 
 
 The strict---or multiplier---topology on $B(G)$ is  taken with respect to the ideal $A(G)$.    In the  following statement, $G$ is amenable, $\psi = j_\beta :A(G) \ra B(H)$, where $\beta:Z\subseteq H \ra G$, is a homomorphism and $\ov{\psi}:B(G) \ra B(H)$ is the unique  extension of $\psi$ to $B(G)$ that is strict-to-weak$^*$ continuous on bounded subsets of $B(G)$, (which exists by \cite[Corollary 5.8]{Il-St}). 
 
 \bp \label{Strict Extension Prop} Let $G$ be amenable, $\psi= j_\beta:A(G) \ra B(H)$ as above,  and let $\vp = j_\alpha:B(G) \ra B(H)$, $\alpha: Y \ra \DBG$, be a homomorphic extension of $\psi$ to $B(G)$. Then $\vp = \ov{\psi}$ if and only if $Y = \ov{Z}$ (and $\alpha{\large |}_Z = \beta$). 
 \ep  
 
 \begin{proof}  We know that $Y \in \Omega(H)$ and by Lemma \ref{Restriction to A(G) Lemma}, $Z = \alpha^{-1}(G) \subseteq Y$ is open and $\beta = \alpha{\large |}_Z$. Therefore $\ov{Z}$ is contained in $Y$. 
 
 Suppose that $\vp = \ov{\psi}$ and assume towards a contradiction that there is an element $y$ belonging to the open set $Y \bs \ov{Z}$. Taking $u \in B(G)$ such that $1 = \l \alpha(y), u\r  = \ov{\psi}(u)(y)$, we can find a compact neighbourhood $K$ of $y$ such that $K$ is contained in $Y \bs \ov{Z}$ and $| \ov{\psi}(u)(x)| \geq 1/2$ on $K$. Define $f \in L^1(H) \subseteq C^*(H)$ by putting $f(x) = 1/\ov{\psi}(u)(x) $ on $K$, 0 otherwise. Letting $(e_\delta)$ be a bounded approximate identity for $A(G)$, $\l \ov{\psi}(u), f \r = \lim_\delta \l \psi(u e_\delta), f\r$ by \cite[Theorem 5.6]{Il-St}. However, 
 $$\l \ov{\psi}(u), f \r = \int_H \ov{\psi}(u)(x) f(x) \, dx = \int_K 1 \, dx >0$$ since $K$ is open and nonempty, and 
 $$\lim_\delta \l \psi(u e_\delta), f\r  = \lim_\delta \int_K \psi(u e_\delta)(x) f(x) \, dx = 0$$ since $\psi(ue_\delta) = j_\beta(ue_\delta)$ vanishes on the subset $K$ of $H \bs Z$. We conclude  that  $Y = \ov{Z}$. 
 
 Conversely, suppose that $\ov{Z} = Y$.  From the above, we know that $\ov{\psi} = j_{\ov{\beta}}$ where $\ov{\beta}:\ov{Z} \ra \DBG$ is a continuous extension of $\beta$ to $\ov{Z}$. But from the first paragraph of this proof, $\alpha$ is also a continuous extension of $\beta$ to $Y = \ov{Z}$, so we must have $\alpha= \ov{\beta}$. Hence, $\vp = j_\alpha = j_{\ov{\beta}} = \ov{\psi}$. 
 \end{proof}

\bc  \label{Strict Extension Cor1}  Let $G$ be amenable.  \bi \item[(i)]  Let $\psi = j_\beta:A(G) \ra B(H)$, $\beta:Z \ra G$, be a homomorphism. Then  $Z$ is open,  $Z \in \Omega(H_d)$, $\ov{Z} \in \Omega(H)$ and $\beta$ extends to a continuous function $\ov{\beta}: \ov{Z} \ra \DBG$ such that $\ov{\beta}(\ov{Z} \bs Z) \subseteq \DBG \bs G$.
\item[(ii)] Let $\vp= j_\alpha:A \ra B(H)$ be a homomorphism. Then $Y_a = \alpha^{-1}(G)$  is closed if and only if $Y_a \in \Omega(H)$, (equivalently, $\vp$ has a Lebesgue decomposition). 

\ei  \ec  

\begin{proof} (i)  By Proposition \ref{Strict Extension Prop}, $\beta: Z \ra G$ has a continuous extension $\ov{\beta}: \ov{Z}\ra \DBG$ such that $\ov{\psi} = j_{\ov{\beta}}: B(G) \ra B(H)$ is the homomorphic extension of $\psi$ that is strict-to-weak$^*$ continuous on bounded subsets of $B(G)$.  By Proposition \ref{FSAlgHomGeneralProp}, $\ov{Z} \in \Omega(H)$. Since $\ov{\psi}_F = \psi = j_\beta$, $Z = \ov{\beta}^{-1}(G)$ by Lemma \ref{Restriction to A(G) Lemma}, so $Z$ is open and $\ov{\beta}$ maps $\ov{Z} \bs Z$ into $\DBG\bs G$.   To see that $Z \in \Omega(H_d)$---where $H_d$ denotes the group $H$ with the discrete topology---just observe that because $B(H) = CB(H) \cap B(H_d)$, $\psi=j_\beta:A(G) \ra B(H_d)$ is also a homomorphism. Hence, $Z = \ov{Z}^{H_d} \in \Omega(H_d)$.   \\
(ii)  By Lemma \ref{Restriction to A(G) Lemma}, $\vp_F = j_{\alpha_a}$ where $\vp_F = \vp{\big |}_{A(G)}$ and $\alpha_a = \alpha{\big |}_{Y_a}$. By part (i). $\ov{Y_a}$ is in $\Omega(H)$, so $Y_a$ is in $\Omega(H)$ when it is closed. The converse is trivial. 
\end{proof}  

 Unfortunately, we do not know if the set  $Z$ from Corollary \ref{Strict Extension Cor1}(i) must always belong to $\Omega(H)$ (equivalently, whether $Z$ is closed) or not. Observe that if $H=\mathbb{R}$, $Z = \mathbb{R} \bs \{0\}$ is open, $Z \in \Omega(H_d)$, $\ov{Z} = H \in \Omega(H)$, but $Z \notin \Omega(H)$ since it is not closed.  So, the conclusions of Corollary \ref{Strict Extension Cor1}(i) do not alone imply that $Z$ is closed. 
 
 \bc \label{Strict Extension Cor2}  Let $G$ be amenable,  $\vp = j_\alpha:B(G) \ra B(H)$ a homomorphism. 
If $\vp$ is spatial, then it is  strict-to-weak$^*$ continuous on bounded subsets of $B(G)$; the converse holds when $\vp$ has a Lebesgue decomposition.  
\ec  

\begin{proof}  Suppose that $\vp$ is spatial. By Proposition \ref{FSAlgHomGeneralProp}, $Y \in \Omega(H)$ and, letting $\vp_F = \vp{\big |}_{A(G)}$,  $\vp_F = j_\alpha$  by Lemma \ref{Restriction to A(G) Lemma}.  Since $Y$ is closed, $\ov{\vp_F} = \vp$ by Proposition \ref{Strict Extension Prop}.  Conversely,   if $\vp$ is strict-to-weak$^*$ continuous and has Lebesgue decomposition $\vp= \vp_a + \vp_s$, then by Lemma \ref{Basic Lemma spatial singular disjoint}(ii), both $\vp$ and $\vp_a$ are  extensions of $\vp {\big |}_{A(G)}$ that are strict-to-weak$^*$ continuous on  bounded subsets $B(G)$. The uniqueness of this extension yields  $\vp = \vp_a$, so $\vp$  is spatial.   
\end{proof}

  \bc \label{Lebesgue Decomp Cor} Let $\vp: A \ra B(H)$ be a homomorphism. In the following situations, $\vp$ has a Lebesgue decomposition $\vp = \vp_a + \vp_s$, written as  $\vp_a = j_{\alpha_a}$, $\vp_s = j_{\alpha_s}$: \bi 
  \item[(i)] $G$ is amenable and $\vp$ is completely bounded. In this case, $\alpha_a :Y_a \ra G$ is a continuous piecewise affine  map on $Y_a \in \Omega(H)$.   
   \item[(ii)] $\vp$ is positive (contractive). In this case,  $\alpha_a :Y_a \ra G$ is a continuous homomorphic or anti-homomorphic (affine or anti-affine) map---a homomorphism (affine map) if $\vp$ is completely positive (completely contractive)---on $Y_a$, an open subgroup (coset) of $H$.   
    \item[(iii)]  $G$ is amenable and $H$ is discrete.  
  \ei  
When $A = B(G)$ and $G$ is amenable, in each case $\vp$ is spatial---equivalently $\vp=\vp_a$
---if and only if $\vp$ is strict-to-weak$^*$ continuous on bounded subsets of $B(G)$.   \ec 
  
  \begin{proof}  Let $Y_a = \alpha^{-1}(G)$, $\alpha_a = \alpha |_{Y_a}: Y_a \ra G$. Then $\vp_F = j_{\alpha_a}$ where $\vp_F$ is the restriction of $\vp$ to $A(G)$. By Proposition \ref{Lebesgue decomposition Prop},  it suffices to show that $Y_a \in \Omega(H)$ in each case:  
  (i) By \cite[Theorem 3.7]{Il-Sp}, $Y_a \in \Omega(H)$, so $\vp$ has Lebesgue decomposition with $\vp_a = j_{\alpha_a}$. That $\alpha_a$ is piecewise affine also follows from \cite[Theorem 3.7]{Il-Sp}.   (ii) If $\vp$ is positive (contractive), $Y_a$ is an open---and therefore closed---subgroup (coset) of $H$ by 
 \cite[Lemma 4.2 and Theorem 5.1]{Pham}, and $\alpha_a$ is a homomorphism or anti-homomorphism by \cite[Theorems 4.3 and  5.1]{Pham}. When $\vp$ is completely positive/completely contractive, this follows from either Proposition \ref{PhamProp5.8} and Remark \ref{PhamProp 5.8 Remark}.2 or \cite[Proposition 5.8]{Pham}. (iii) In this case, $Y_a$ is trivially closed and therefore $Y_a \in \Omega(H)$ by Corollary \ref{Strict Extension Cor1}.
  The last statement follows from Corollary \ref{Strict Extension Cor2}. 
\end{proof}

 Thus, when $\vp: A \ra B(H)$ is spatial, $\vp$ is completely described in cases (i) and (ii) above. On the other hand, very little is known about non-spatial homomorphisms.  It may be natural to wonder if $\vp_a$ and $\vp_s$ can both be nontrivial in the Lebesgue decomposition of a homomorphism $\vp$, but examples of such homomorphisms are given by the proofs of \cite[Theorem 4.6.1(b)]{Rud} and 
\cite[Proposition 5.1(ii)]{Il-St}.  In the next section, we will show how to build a wide class of such examples. (For now, just observe that if $\vp_a, \vp_s: A \ra B(H)$ are nontrivial disjoint homomorphisms that are, respectively, spatial and singular,  then $\vp = \vp_a + \vp_s$ is a homomorphism by Lemma \ref{Disjoint Homom Lemma}.) We will show that even a completely positive homomorphism can have nontrivial spatial and singular parts in Section 4.    
 
  \br \rm  Let $G$ be amenable, $\vp = j_\alpha: B(G) \ra B(H)$ any homomorphism.  Then $\vp$ always has a unique decomposition 
  \beq \label{Strict Decomp Eqn} \vp = \vp_1 + \vp_2
  \eeq
  where $\vp_1$ is a homomorphism that is strict-to-weak$^*$ continuous on bounded subsets of $B(G)$ and $\vp_2$ is a singular homomorphism that is disjoint from $\vp_1$.  Moreover, if $Y_a = \alpha^{-1}(G)$ and $\alpha_a = \alpha{\big |}_{Y_a}$, $\vp_1 = \ov{\vp_F} = j_{\ov{\alpha_a}}$ where $\ov{\alpha_a}:\ov{Y_a} \ra \DBG$ is a continuous extension of $\alpha_a$ to $\ov{Y_a}\in \Omega(H)$.  (This follows by replacing $\vp_F= j_{\alpha_a}$ and $Y_a$ in the proof of (b) implies (a) of Proposition \ref{Lebesgue decomposition Prop} with $\ov{\vp_F} =  j_{\ov{\alpha_a}}$ and $\ov{Y_a}$, which is possible by Proposition \ref{Strict Extension Prop}.) Observe that the Lebesgue decomposition of $\vp$ agrees with (\ref{Strict Decomp Eqn}) when it exists.   \er

 \section{Construction of  Fourier--Stieltjes algebra homomorphisms}

 We begin with the case of completely positive homomorphisms. 
 
 \smallskip  
 
 Let $p \in {\cal L}(\H)$ be a non-zero projection, 
  $$\bB_p^\sharp = \{ T \in {\cal L}(\H):  Tp = T \ {\rm and} \  T^\sharp p = T^\sharp, \ T T^\sharp = T^\sharp T = p \ {\rm for  \ some} \ T^\sharp \in {\cal L}(\H) \}$$ and $$\bB_p = \{ T \in {\cal L}(\H):  T T^* = T^* T = p\}.$$ 
 One can readily show that $\bB_p^\sharp$ is a group with identity $p$ and for $T \in \bB_p^\sharp$,  its unique inverse in   $\bB_p^\sharp$ is $T^\sharp$. 

 Let $\K = p \H$, $p_\K: \H \ra \K: \xi \mapsto p \xi$, $E_\K: \K \hookrightarrow \H$ the embedding map. Then 
$T \mapsto p_\K T E_\K \left(= E_\K^* T E_\K = p\left(T\big{|}_\K\right)\right)$ defines 
  a normal unit-preserving completely positive mapping of ${\cal L}(\H)$ onto ${\cal L}(\K)$ whose restriction to the von Neumann subalgebra $p {\cal L}(\H)p$ of ${\cal L}(\H)$, 
 $$\Phi: p {\cal L}(\H) p \ra {\cal L}(\K): T\mapsto p_\K T E_\K,$$ is a (completely isometric, normal) $*$-isomorphism onto ${\cal L}(\K)$ with inverse 
 $$\Phi^{-1}:  {\cal L}(\K)  \ra p {\cal L}(\H)p : S\mapsto E_\K S p_\K.$$   We use $GL(\K)$ to denote the semitopological group of invertible operators on $\K$  with the relative $\sigma = \sigma({\cal L}(\K), {\cal L}(\K)_*)$-topology, and  ${\cal U}(\K)$ to denote the topological group of unitary operators on $\K$ with the SOT-topology. Recall that the $\sigma$-, WOT- and SOT-topologies all agree on ${\cal U}(\K)$. 
 
 \blem \label{BpLem} The set $\bB_p^\sharp $ with the relative $\sigma$-topology from ${\cal L}(\H)$ is a semitopological group  such that 
 $$\Phi: \bB_p^\sharp \ra GL(\K): T \mapsto p_\K T E_\K$$ is a norm-preserving topological isomorphism of $\bB_p^\sharp$ onto $GL(\K)$. Also, $\bB_p$ is a topological subgroup of $\bB_p^\sharp$ that is topologically isomorphic with ${\cal U}(\K)$ via the isomorphism $\Phi$. 
   \elem 
   
   \begin{proof}  Since multiplication in ${\cal L}(\H)$ is separately $\sigma$-continuous, $\bB_p^\sharp$ is a semitopological group and, since $\bB_p^\sharp$ is contained in $p{\cal L}(\H)p$,   $\Phi$ is  multiplicative  and norm preserving on $\bB_p^\sharp$ by the remarks preceding the lemma. Moreover, for $T \in \bB_p^\sharp$, $\Phi(T) \Phi(T^\sharp) = \Phi(TT^\sharp) = \Phi(p) = {\rm id}_\K$ and $\Phi(T^\sharp) \Phi(T)  = {\rm id}_\K$, so $\Phi$ is a continuous homomorphism of $\bB_p^\sharp$ into $GL(\K)$.  Since $\Phi^{-1}(S)  = E_\K S p_\K$ clearly maps  $GL(\K)$ continuously into $\bB_p^\sharp$, we obtain the first statement.  
 Suppose now that $T$ is in  $ \bB_p$. Then  $T$ has polar decomposition $T = Up$, so $Tp= T$; therefore $pT= TT^*T = Tp=T$, whence  $T^*p  = (pT)^* =  T^*$.  Hence, $T \in \bB_p^\sharp$ (and $T^\sharp = T^*$).  For $T \in \bB_p$ and $S \in {\cal U}(\K)$, 
  $$\Phi(T)^* = \Phi(T^*) = \Phi(T^\sharp) = \Phi(T)^{-1},  \  \Phi^{-1}(S)^* = \Phi^{-1}(S^*) = \Phi^{-1}(S^{-1}) = \Phi^{-1}(S)^\sharp,$$ 
  so $\Phi$ maps $\bB_p$ into ${\cal U}(\K)$ and $\Phi^{-1}$ maps ${\cal U}(\K)$ into $\bB_p$. This completes the proof. 
   \end{proof}

 \br  \rm  If we represent operators $T \in {\cal L}(\H) = {\cal L}(\K \oplus \K^\perp)$ as $2\times 2$ matrices $ T = [T_{ij}]$,    Lemma  \ref{BpLem} says that 
$$\bB_p^\sharp = \left\{  \begin{bmatrix}
    S       & 0  \\
   0       &0
\end{bmatrix} : S \in GL(\K)       \right\} \quad {\rm and}  \quad   \bB_p = \left\{  \begin{bmatrix}
    U       & 0  \\
   0       &0
\end{bmatrix} : U \in {\cal U}(\K)       \right\}. $$

\er

 \bi \item[] \it   Throughout the remainder of this section,   $A=A_\pi$ is a closed translation-invariant subalgebra of $B(G)$  for which $\Delta(A)$ is a $*$-semigroup comprised of norm-one elements (e.g., when $A$ is unital or $A= A(G)$).  \rm
\ei 

 Let  $p \in \Delta(A) \subseteq VN_\pi \subseteq {\cal L}(\H_{\pi})$ be a projection, 
 $${\bD}_p = \{x \in \Delta(A): x^*x = x x^* = p\}.$$ 
Then $\bD_p = \bB_p \cap \Delta(A)$, so $\bD_p$ is a topological group by Lemma \ref{BpLem} since $\Delta(A)$ is a $*$-semigroup.  

\br \rm \label{Remark zF}  When $A=B(G)$,  this is found in Walter's paper \cite{Wal2} where the notation $G_p$ is employed.   Walter, ibid.,  showed   $\Delta(B(G))$ always contains a central idempotent $z_F$ such that  $z_F = \min \{ s \in \DBG:  s \geq 0\}$ and  $z_F \neq e_G$ when $G$ is noncompact.  (Corollary \ref{Walter min idempotent etc Corollary}, below, provides simple new proofs of these and other results about $z_F$ from \cite{Wal2}.)        \er 

 Observe that by Corollary \ref{Wal1Cor},     if $A$ contains $A(G)$, then  $\bD_p \subseteq \Delta(A) \bs G$ when $p \neq e_G$; otherwise $\bD_p= G$. A main point of the following proposition is that inversion in any subgroup of $\Delta(A)$ is given by involution.

 \bp \label{SpectralSubgroupsProp}  Let $D$ be a subgroup of $\Delta(A)$ with identity $p$. Then $p$ is a projection and $D$ is a subgroup of $\bD_p$.  Thus, for each $x \in D$, 
 $x^*x = xx^*=p$, $xp = px = x$, and $D$ is a topological group that is topologically isomorphic to a subgroup of ${\cal U}(\K)$ where  $\K = p \H_\pi$.  \ep
 
 \begin{proof} As a norm-one idempotent, $p$ is a projection in ${\cal L}(\H_\pi)$ and it is clear that  $D$ is  a subgroup of $\bB_p^\sharp$. Let $x \in D$. Since $x, x^\sharp \in \Delta(A)$ and $\Phi: \bB_p^\sharp \ra GL(\K)$ is norm-preserving, $\Phi(x), \Phi(x^\sharp) = \Phi(x)^{-1} \in GL(\K)_{\| \c \| = 1}$, whence $\Phi(x)\in {\cal U}(\K)$. Therefore $x = \Phi^{-1}(\Phi(x)) \in \bB_p \cap \Delta(A) = \bD_p$ by Lemma \ref{BpLem}.  Uniqueness of the inverse  $x^\sharp$ in $\bB_p^\sharp$  gives $x^\sharp = x^*$, (i.e., the inverse of $x$ in $D$ is $x^*$).    \end{proof}

 Let $\alpha : H \ra \Delta(A)$ be a continuous homomorphism, ($\alpha(st) = \alpha(s)\alpha(t)$, $s, t \in H$). Then $\alpha(H)$ is a subgroup of $\Delta(A)$ with identity $p = \alpha(e_H)$ so---by Proposition \ref{SpectralSubgroupsProp}---$\alpha(H)$ is a subgroup of   $\bD_p$.   Letting $\K= p \H_\pi$, by Lemma \ref{BpLem}  and Proposition \ref{SpectralSubgroupsProp}, $$\Phi: \bD_p \ra {\cal U}(\K): x \mapsto p_\K x E_\K$$ is a continuous $*$-monomorphism, 
  so $$\pi_\alpha \coloneqq \Phi \circ \alpha : H \ra {\cal U}(\K)$$ is a continuous unitary representation of $H$ on $\K$.  When $A= B(G)= A_{\omega_G}$, we write $\omega_\alpha$ in place of $(\omega_G)_\alpha$.
 
 \bp \label{CPHomProp} Let $\alpha: H \ra \Delta(A)$ be a continuous homomorphism with $\alpha(e_H) = p$. Then the map $j_\alpha: A \ra B(H)$ given by $j_\alpha u(h) = \l \alpha(h), u\r_{VN_\pi-A_\pi}$ is a well-defined completely positive homomorphism such that for each $u = \xi*_\pi \eta$, $j_\alpha u = p_\K \xi *_{\pi_\alpha} p_\K \eta$ and $\ja^*:W^*(H) \ra VN_\pi$ is a normal $*$-homomorphism extending $\alpha$. Moreover: 
 \bi 
 \item[(i)] $\ja$ is a complete quotient mapping of $A$ onto the (completely contractive)  Banach subalgebra $A_{\pi_\alpha}$ of $B(H)$ such that for each $v \in A_{\pi_\alpha}$ there is a $u \in A$ such that $\ja (u) = v$ and $\| u \| = \| v\|$; 
 \item[(ii)]  $A_{\pi_\alpha}$ separates points of $H$ when  $\alpha$ is one-to-one; and
 \item[(iii)] when $A=B(G)$,  for each $v \in A_{\omega_\alpha}$ there exist $\xi, \eta \in \K$ such that $v =  \xi *_{\omega_\alpha} \eta$ and $\| v \| = \| \xi \| \| \eta\|$. 
    \ei 
 \ep   
 
 \begin{proof}  Let $u = \sum \xi_n *_\pi \eta_n \in A$ with $\|u \| = \sum\|\xi_n\|\|\eta_n\|$ \cite[Theorem 2.2]{Ars}. Then $\sum \|p_\K \xi_n\|\|p_K\eta_n\| \leq \|u\|$, so $v= \sum p_\K \xi_n *_{\pi_\alpha}  p_\K \eta_n \in A_{\pi_\alpha} \subseteq B(H)$ with $\| v\| \leq \|u\|$  \cite[Theorem 2.2]{Ars}, and for $h \in H$,   
 \beqs j_\alpha u(h) & = & \l \alpha(h), \sum \xi_n*_\pi \eta_n\r  =  \sum \l\alpha(h) \xi_n | \eta_n \r \\
 & = & \sum  \l p \alpha (h) p \, \xi_n | \eta_n \r   = \sum  \l  E_\K p_\K \alpha(h) E_\K p_\K  \xi_n | \eta_n \r\\
 &  = &\sum  \l \Phi(\alpha(h)) p_\K \xi_n | E_\K^* \eta_n \r  = \sum  p_\K \xi_n *_{\pi_\alpha} p_\K \eta_n(h) = v(h). 
 \eeqs 
 Thus, $\vp = j_\alpha$ is a contraction mapping $A=A_\pi$ into $A_{\pi_\alpha}$ and, since $\alpha(H)$ is contained in $\Delta(A)$, $\vp $ is a homomorphism of $A$ into $B(H)$. Observe that $\vp^*: W^*(H) \ra VN_\pi$ satisfies $\vp^*(\omega_H(h)) = \alpha(h)$, so $\vp^*$ is a $*$-homomorphism on $\l \omega_H(H)\r$, which is $\sigma$-dense in $W^*(H)$. By separate $\sigma$-continuity of multiplication in $W^*(H)$ and $\sigma$-continuity of $\vp^*$ and the involution, $\vp^*$ is a $*$-homomorphism. Hence, $\vp$ is  completely positive.  
 
 Supposing now that $v = \sum \xi_n*_{\pi_\alpha}\eta_n \in A_{\pi_\alpha}$ with  $\|v \| = \sum \| \xi_n\| \|\eta_n\|$, $u = \sum E_\K \xi_n*_{\pi} E_\K \eta_n \in A$ with   $\|u \| \leq \sum \| E_\K \xi_n\| \|E_K \eta_n\| = \|v\|$  and $\ja u = v$; since $\ja$ is contractive  $\|v\| = \| \ja u \| = \| u \|$.  Hence,  $\vp_0: A \ra A_{\pi_\alpha}: u \ra \ja(u)$ is a quotient mapping and, since  
  $\vp_0^* : VN_{\pi_\alpha} \ra VN_\pi$ maps $\pi_\alpha(h)$ to $\alpha(h)$, it, like $\vp^*:W^*(H) \ra VN_\pi$ above, is   a $*$-homomorphism on  $VN_{\pi_\alpha}$;  hence $\vp_0^*$ is, in fact, a (completely isometric) $*$-isomorphism and therefore $\vp_0$ is a complete quotient mapping.  This proves (i).  
  
   When $A=B(G)$, we  can  take $\xi, \eta \in \H_{\omega_G}$ such that $u = \xi *_{\omega_G}\eta$ and $\|u \| = \| \xi \| \|\eta\|$ (this follows, e.g.,  from \cite[(2.14)]{Eym}).  Then $v = \ja u =p_\K \xi *_{\omega_\alpha} p_\K \eta$ and 
 $\| v \| \leq \| p_\K \xi \| \| p_\K \eta\| \leq \| \xi \| \| \eta \| = \|u \| = \| v \|$, so $\| v\| =\| p_\K \xi \| \| p_\K \eta\|$. This is (iii).

  Finally, since $\pi_\alpha$ is faithful when $\alpha$ is one-to-one, statement (ii) is obvious.   \end{proof} 
  
  \bc \label{CPHomCor} If $H_0$ is an open subgroup of $H$ and $\alpha:H_0 \ra \Delta(A)$ is a continuous homomorphism, then  $\ja: A \ra B(H)$ is a completely positive homomorphism. \ec 
 \begin{proof}  In the second paragraph of the proof of \cite[Proposition 3.1]{Il-Sp}, the authors showed that the map $u \mapsto u^\circ$ is a completely positive homomorphism of $B(H_0)$ into $B(H)$, where $u^\circ(h) = u(h)$ for $h \in H_0$ and $u^\circ(h) = 0$ for $h \in H\bs H_0$. By Proposition \ref{CPHomProp}, $\ja$ is hence a composition of completely positive homomorphisms and is therefore, itself, completely positive. 
 \end{proof} 
 
  \br  \label{CPHomProp Remarks} \rm    1. If $\alpha: H \ra \Delta(A)$ is an anti-homomorphism, then, using the fact that  involution is a positive operator, one can use Corollary \ref{CPHomCor} to show that $\vp= \ja$ is a positive (though not necessarily completely positive) homomorphism.

 \noindent 2. If $p=e_G$, the continuous homomorphism $\alpha$ maps $H$ into  $G$, and it was already known in this case that $j_\alpha$ is a completely positive homomorphism of $B(G)$ into $B(H)$ \cite[Corollary 3.2]{Il-Sp}. 
 
 \noindent 3.  When $A = A(G)$, the converse of Corollary \ref{CPHomCor} holds. 
  In Section 4,   we will  exhibit a large collection of completely positive homomorphisms that do not take this form.  \rm 
 \er 

 
 
 
 
When $\alpha: H \ra G$ is a continuous homomorphism and $A_\tau$ is a closed translation-invariant subspace of $B(G)$, \cite[Proposition 2.10]{Ars} describes the image of $A_\tau$ under $\ja$ (also see Sections 2.2 and 2.8 of \cite{Kan-Lau}).  We now  extend  \cite[Proposition 2.10]{Ars} (and  refine Proposition \ref{CPHomProp} when $A=B(G)$)  by  describing the image of $A_\tau$ under $\ja$ when $\alpha: H \ra \DBG$ is any homomorphism.    Part (iii) of Proposition \ref{Arsac2.10GenProp}  will be used in  Section 5. 

 Let $\{\tau, \H_\tau\}$ be a continuous unitary representation of $G$.  Letting $\iota: A_\tau \hookrightarrow B(G)$, $\tau_\omega = \iota^*: W^*(G) \ra VN_\tau$ is the normal $*$-representation of $W^*(G)$ onto $VN_\tau$ that extends $\tau$ on $G$.  We write $\dot{\xi}$ when viewing an element $\xi$ in $\H_\tau$ as an element of $\ov{\H_\tau}$, the conjugate Hilbert space of $\H_\tau$,  and for $B \in {\cal L}(\H_\tau)$ we define $\dot{B} \in {\cal L}(\ov{\H_\tau})$ by $\dot{B}(\dot{\xi}) = (B \xi)^{\dot{}}$. The conjugate representation $\{\ov{\tau}, \ov{\H_\tau}\}$ of $\{ \tau, \H_\tau \}$ can then be defined as $\ov{\tau}(s) = (\tau(s))^{\dot{}}$ $(s \in G)$. Letting $G^{\cal E}$ denote the closure of $G$ in $\DBG$---the Eberlein compactification of $G$, e.g., see \cite{Spr-Sto}---we observe that 
 \beq \label{ConjugateRepFormula}   (\ov{\tau})_\omega(x) = (\tau_\omega(x))^{\dot{}} \quad {\rm whenever} \quad x \in G^{\cal E}. \eeq
 This follows because $B \mapsto \dot{B}$  is $\sigma$-continuous and the maps $\tau_\omega$ and $(\ov{\tau})_\omega$ are normal mappings of $W^*(G)$ into ${\cal L}(\H_\tau)$ and ${\cal L}(\ov{\H_\tau})$, respectively.

  Let  $\alpha: H \ra \DBG$ be a continuous homomorphism with $p = \alpha(e_H)$. By Proposition \ref{CPHomProp}, $\ja:B(G) \ra B(H)$ is a completely positive homomorphism.  Let $q = \tau_\omega(p)$ and assume that $q$ is nonzero.    Then $\tau_\omega \circ \alpha (H) \subseteq \bB_q \subseteq {\cal L}(\H_\tau)$, so $$\tau_\alpha \coloneqq \Phi_q \circ \tau_\omega \circ \alpha: H \ra {\cal U}(\K)$$ is a continuous unitary representation of $H$, where we are now taking $\K = q \H_\tau$ and $\Phi_q: \bB_q \ra {\cal U}(\K): x \ra p_\K x E_\K$.  
 
 \bp \label{Arsac2.10GenProp} Let $\alpha: H \ra \DBG$ be a continuous homomorphism with $p = \alpha(e_H)$, $\ja: B(G) \ra B(H)$ the associated homomorphism of Fourier-Stieltjes algebras. 
 \bi  \item[(i)] If $q = \tau_\omega(p) \neq 0$, then $\ja$ maps $A_\tau$  as a completely positive complete quotient mapping onto $A_{\tau_\alpha}$, and for each $v$  in $A_{\tau_\alpha}$ there is a $u$ in $A_\tau$ such that $\ja u =v $ and $\| u \| = \| v\|$; if $\tau_\omega(p) = 0$, then $\ja$ maps $A_\tau$ to $\{0\}$. 
 \item[(ii)] If $\ja$ is weak$^*$-continuous,  then $\ja$ maps $B_\tau$, the weak$^*$-closure of $A_\tau$ in $B(G)$, onto  $B_{\tau_\alpha}$. 
 \item[(iii)]  If $\alpha(H)$ is contained in $G^{\cal E}$, then $\ja(\ov{u}) = \ov{\ja(u)}$ for every $u \in B(G)$. 
 \ei  
  \ep   
 
 \begin{proof}   Observe that for  $u = \xi *_\tau \eta$ in $A_\tau$ and $h \in H$, 
 \beq  \label{ArsacPropEqn} \ja u(h) = \l \alpha(h), \xi*_\tau \eta \r_{W^*-B(G)} = \l \tau_\omega(\alpha(h)), \xi*_\tau \eta \r_{VN_\tau-A_\tau}  = \l \tau_\omega(\alpha(h)) \xi | \eta\r_{\H_\tau}.    \eeq
 By modifying the calculations found in the proof of Proposition \ref{CPHomProp}, it follows   that  $\ja (u) = p_\K \xi*_{\tau_\alpha} p_\K \eta$ when $\tau_\omega(p) \neq 0$. The remainder of the proof of (i) is essentially the same as that of Proposition \ref{CPHomProp}, and the proof of  \cite[Lemma 4.1]{Il-St}  now yields (ii).  To establish (iii), let $u \in B(G)$, say $u = \xi *_\tau \eta$ for a continuous unitary representation $\{\tau, \H_\tau\}$ of $G$, and let $h \in H$.  Observe that $\ov{u} = \dot{\xi} *_{\ov{\tau}} \dot{\eta}$ and we are assuming that $\alpha(h) \in G^{\cal E}$, so (\ref{ConjugateRepFormula}) and (\ref{ArsacPropEqn}) give 
 \beqs \ja(\ov{u})(h) & = & \ja(\dot{\xi} *_{\ov{\tau}} \dot{\eta})(h)  = \l(\ov{\tau})_\omega(\alpha(h)) \dot{\xi} | \dot{\eta}\r \\
 & = & \l  \tau_\omega(\alpha(h))^{\dot{}} ( \dot{\xi}) | \dot{\eta}\r  = \ov{\l \tau_\omega(\alpha(h)) \xi | \eta\r } \\
 & = &  \ov{\ja u(h)}
 \eeqs
  as needed.  \end{proof}

\br \rm 1. If $p = \alpha(e_H)=e_G$, then $\alpha(H) \subseteq G$ and $q = \tau_\omega(e_G) = {\rm id}_{\H_\tau}$, so $\bB_q = {\cal U}(\H_\tau)= {\cal U}(\K)$ and $\tau_\alpha (h) = \Phi_q \circ \tau_\omega \circ \alpha(h) = \tau \circ \alpha(h)$. Hence, \cite[Proposition 2.10]{Ars} is a special case of Proposition \ref{Arsac2.10GenProp}.

\noindent 2.  When $\tau = \omega_G$, $\tau_\alpha = \omega_\alpha$ (with no ambiguity in the notation), so in the case that $A=B(G)$, parts of  Proposition \ref{CPHomProp} are also contained in part (i) of Proposition \ref{Arsac2.10GenProp}.

\noindent 3. Replacing $A=B(G)$ in Proposition \ref{Arsac2.10GenProp} with any algebra $A=A_\pi$ satisfying the general assumptions of this section,  and taking $A_\tau \subseteq A_\pi$, one can also prove a more general version of Proposition \ref{CPHomProp}(i).  However, this is not needed in the sequel.

\noindent 4. In the special case that the homomorphism $\alpha$ maps $H$ into $G$, it was shown in \cite{Il-St}  that $\ja$ is weak$^*$-continuous exactly when $\alpha$ is an open map. In general, we do not yet know when $\ja$ is weak$^*$-continuous. 
  \er


We now turn our attention to the construction of  completely contractive homomorphisms and, then,  completely bounded  homomorphisms.    

\medskip

 As a von Neumann algebra, $VN_\pi$ is a completely contractive dual Banach algebra, and therefore $A= A_\pi$ is a completely contractive $VN_\pi$-submodule of $VN_\pi^*$ via  $$\l y, u \c x \r = \l xy, u \r \ \ {\rm and}  \ \ \l y, x\c u \r = \l yx, u \r \ \ {\rm for } \ x,y \in VN_\pi, \ u \in A.$$
 For $x \in VN_\pi$, let $$\ell_x, r_x : A \ra A: u \mapsto  \ell_x u= u \c x, r_x u = x \c u.$$
 Then $\ell_x, \, r_x$ are completely bounded with $\| \ell_x\|_{cb}=  \| r_x \|_{cb} = \|x \|$ and for any $x, y$, $\ell_x \circ r_y = r_y \circ \ell_x$.  We will need the following lemmas. 
  
 \blem \label{TranslationCCHomLem}  For $x \in \Delta(A)$, $\ell_x, r_x : A\ra A$ are completely contractive homomorphisms. If we view $\ell_x, \, r_x$ as homomorphisms of $A$  into $B(G)$, $\ell_x = j_{\lambda_x}$, $r_x = j_{\rho_x}$ where $\lambda_x, \, \rho_x : G \ra \Delta(A): g \mapsto x\pi(g), \, \pi(g) x.$ 
 \elem  
 
 \begin{proof} Since $\| x\|=1$, $\ell_x$ is a complete contraction. As $\ell_x u(g)  = u \c x (g) = \l x\pi(g), u \r$
 and $x \pi(g) \in \Delta(A)$, $\ell_x$ is a homomorphism.   If we view $\ell_x $ as  a  homomorphism into $B(G)$, 
 $\ell_x^*(g) = x \pi(g) \neq 0$, so by Proposition \ref{FSAlgHomGeneralProp}, $\ell_x =  j_{\lambda_x}$. 
 The same argument gives the statement for  $r_x$.  
 \end{proof} 
 
  Observe that  $\ell_{\pi(x)} u(g)  = u \c x(g) = u (xg)$ when $x, g \in G$, which is consistent with the usual notation. 
We remark that although we have not found it helpful to use the Gelfand transform to identify $B(G)$ with its algebraic copy $\{ \widehat{u}: u \in B(G)\}$ in $CB(\DBG)$, we could write $\ell_x u(y) = u \c x(y) = u (xy)$ for $u \in B(G)$ and $x,y \in \DBG$ if we did make this identification.

 \blem \label{CCHomLem}  Let   $\vp = \ja: A \ra B(H)$ be a (completely contractive) homomorphism where $\alpha: Y \ra \Delta(A)$, and for $y_0, y_1 \in H$ and $x_0, x_1  \in \Delta(A)$, define $\beta: y_0^{-1}Yy_1^{-1} \ra \Delta(A)$  by $\beta(h) = x_0 \alpha(y_0hy_1)x_1$.  Then $j_\beta: A \ra B(H)$ is a  well-defined (completely contractive) homomorphism that factors as $j_\beta = r_{y_1} \circ \ell_{y_0}\circ \ja \circ \ell_{x_0}\circ r_{x_1}$.
 \elem
 
 \begin{proof}  By Lemma \ref{TranslationCCHomLem}, $\psi =  r_{y_1} \circ \ell_{y_0}\circ \ja \circ \ell_{x_0} \circ r_{x_1} : A \ra B(H)$ is a (completely contractive) homomorphism. For $h \in H$, 
 \beqs \psi^*(h) & = & r_{x_1}^* \circ  \ell_{x_0}^* \circ \ja^* \circ \ell_{y_0}^* \circ r_{y_1}^*(h)  = x_0 \ja^*(y_0hy_1)x_1 \\
 & = &   \displaystyle{\left \{ \begin{array}{ll}
                           x_0\alpha(y_0hy_1)x_1  & \mbox{{\rm if} $y_0 h y_1 \in Y$}\\
                          0 & \mbox{{\rm if} $ y_0 h y_1  \in H\bs Y $}
                            \end{array}
               \right .  = \left \{ \begin{array}{ll}
                           \beta(h)  & \mbox{{\rm if} $h \in y_0^{-1}Yy_1^{-1}$}\\
                          0 & \mbox{{\rm if} $y \in H \bs y_0^{-1}Yy_1^{-1}$},
                            \end{array}
                 \right. }      
 \eeqs        
 so $\{h \in H: \psi^*(h) \neq 0 \} = y_0^{-1}Yy_1^{-1}$ and $\psi^*\large{|}_{y_0^{-1}Yy_1^{-1}} = \beta$. Hence, $\psi = j_\beta$ by Proposition \ref{FSAlgHomGeneralProp}.
                 \end{proof}

If $\alpha$ maps a coset $Y\subseteq H$ into a subgroup $F$ of a semigroup $S$, we will say that $\alpha$ is \it affine  \rm if it is affine when viewed  as a mapping of $Y$ into $F$. 
                 
\bp \label{CCProp} Let $\alpha: Y \ra \Delta(A)$ be a continuous affine map defined on an open coset $Y$ of $H$. Then  $\ja: A \ra B(H)$ is a completely contractive homomorphism. 
 \ep

 \begin{proof}  Writing $Y$ as  $Y= y_0H_0$ for some open subgroup $H_0$ of $H$,   there is a continuous homomorphism $\beta:H_0 \ra \Delta(A)$ such that $\alpha(y) = \alpha(y_0) \beta(y_0^{-1}y)$.  By Corollary \ref{CPHomCor}, $j_\beta$ is a complete contraction so,  by Lemma \ref{CCHomLem},  $\ja = \ell_{y_0^{-1}} \circ j_\beta \circ \ell_{\alpha(y_0)}$  is a completely contractive homomorphism.
 \end{proof}


Let $p$ be a projection in $\Delta(A)$, $Y \in \Omega(H)$, and let $\alpha: Y \ra \bD_p$ be a continuous piecewise-affine map. By \cite[Proposition 4.4]{Il}, there are pairwise disjoint sets $Y_1, \ldots, Y_n \in \Omega_0(H)$ such that $Y = \bigcupdot_1^n Y_i$ and for each $i$, $\alpha_i:= \alpha \big{|}_{Y_i}$ has a continuous affine extension $\gamma_i$ mapping $E_i = {\rm Aff}(Y_i)$---the smallest (open) coset in $H$ containing $Y_i$---into $\bD_p$. (We note that local compactness of $\bD_p$ is not required for this description of a piecewise affine map.) By Proposition \ref{CCProp},  each $j_{\gamma_i}: A \ra B(H)$ is a completely contractive homomorphism and, since $1_{Y_i}$ is an idempotent in $B(H)$, $\mathfrak{1}_i(u) = u 1_{Y_i}$ defines a completely bounded automorphism of $B(H)$. Hence, each $j_{\alpha_i} = \mathfrak{1}_i \circ j_{\gamma_i}: A \ra B(H)$ is a completely bounded homomorphism and therefore, by Lemma \ref{Disjoint Homom Lemma}, $j_\alpha = j_{\cupdot \alpha_i} = \sum_1^n j_{\alpha_i}: A \ra B(H)$ is a completely bounded homomorphism. 

We shall say that a map $\alpha: Y \ra \Delta(A)$ \it is $pw^2$-affine  \rm (for piecewise-piecewise-affine) if there are projections $p_1, \ldots, p_n$ in $\Delta(A)$ and piecewise-affine maps $\alpha_i: Y_i \ra \bD_{p_i}$ such that $Y = \bigcupdot_1^n Y_i$ and $\alpha = \cupdot \alpha_i$. Assuming further that  $\alpha$ is continuous, we know from the above discussion and Lemma \ref{Disjoint Homom Lemma} that $\ja = j_{\cupdot \alpha_i}: A \ra B(H)$ is a completely bounded homomorphism. We have established:

\bp  \label{PWAffineCBProp} Let $\alpha: Y \ra \Delta(A)$ be a  pw$^2$-affine map.  Then $\ja: A \ra B(H)$ is a completely bounded homomorphism.  
\ep

\bex \label{Rudin Ex as pw2 affine} \rm Suppose that $G$ is non-compact. It should be clear that Proposition \ref{PWAffineCBProp} makes it easy to construct a large collection of (non-spatial) completely bounded homomorphisms from  $B(G)$ into $B(H)$. It seems that the only previously known examples of non-spatial homomorphisms were provided by the construction found in the proof of \cite[Proposition 5.1(ii)]{Il-St} and, in the abelian case,  Rudin's construction from the proof of \cite[Theorem 4.6.1(b)]{Rud}.  Let us observe that the constructions from  \cite{Il-St} and \cite{Rud}  can be expressed in terms of pw$^2$-affine maps: 

Suppose that $G$ is non-compact and $\alpha_a: Y_a \ra G$ is any continuous piecewise-affine map where $Y_a\in \Omega(H)$ is a proper subset of $H$.  Pick any element $x_s$ in $\bD_{z_F}$ (or in any other subgroup $\bD_p$ of $\Delta(A)$  with $p \neq e_G$) and consider the map $\alpha_s: H \bs Y_a \ra \bD_{z_F}: y \mapsto x_s$. Since we can express $Y_s= H\bs Y_a$ as $Y_s = Y_1 \cupdot \cdots \cupdot Y_n$ with each $Y_i \in \Omega_0(H)$, it is clear that $\alpha_s$ is piecewise affine. Hence, $\alpha:= \alpha_a \cupdot \alpha_s: H \ra G \cupdot \bD_{z_F} \subseteq \DBG$  is a continuous pw$^2$-affine map and therefore, by Proposition \ref{PWAffineCBProp}, $\vp= \ja$ is a non-spatial completely bounded homomorphism. (Observe that since $\vp= \ja$ is a non-spatial extension of $\vp \large{|}_{A(G)}= j_{\alpha_a}: A(G) \ra B(H)$, this also provides a new proof of   \cite[Proposition 5.1(ii)]{Il-St} and \cite[Theorem 4.6.1(b)]{Rud}.)
\eex

\br  \rm 1. Let $F$ be a group. In \cite{Spr}, N. Spronk introduced the notion of  a \it mixed \rm piecewise-affine map $\alpha: Y \ra F$; the definition is the same as that of a piecewise-affine map provided above, except that the maps $\gamma_i$ can either be affine or anti-affine. (Anti-affine maps share the same relationship with anti-homomorphisms that affine maps share with homomorphisms.) Defining mixed pw$^2$-affine maps $\alpha:Y \ra \Delta(A)$ in the obvious way, it follows from Remark  \ref{CPHomProp Remarks}.1 that any such map determines a bounded (but not necessarily completely bounded) homomorphism $\ja: A \ra B(H)$.  (Spronk conjectured \cite[Conjecture 4.8]{Spr} that the mixed piecewise-affine maps $\alpha: Y \ra G$ are in one-to-one correspondence with the collection of all  homomorphisms $\vp=\ja:A(G) \ra B(H)$.)

\noindent 2. It is natural to ask if the converse to Proposition \ref{PWAffineCBProp} holds when $G$ is amenable and $\Delta(A)$ is a Clifford semigroup, (meaning, by Proposition \ref{SpectralSubgroupsProp}, that each element in $\Delta(A)$ belongs to some $\bD_p$). On the other hand, if there exists an element $x_0 \in \Delta(A)$ that does not belong to a subsemigroup of $\Delta(A)$, then the constant map $\beta: H \ra \Delta(A): h \mapsto x_0$ is not, according to our definition, pw$^2$-affine.  However, $j_\beta: A \ra B(H)$ is nonetheless a completely contractive homomorphism, since  $\beta(h) = x_0 \alpha(h)$ where $\alpha :H \ra \pi(G) \subset \Delta(A)$ is the trivial homomorphism (see Corollary  \ref{CPHomCor} and Lemma \ref{CCHomLem}). Unfortunately, we are not aware of a known example of $A$ for which $\Delta(A)$ is not a Clifford semigroup.   
\er

\section{More on the structure of $\Delta(A)$ and  applications}


 Let $AP(G)$ denote the algebra of continuous almost periodic functions on $G$, $(G^{ap}, \delta^{ap})$ the almost periodic compactification of $G$---also known as the Bohr compactification of $G$---and let $A_\F = A_\F(G) = B(G) \cap AP(G)$. 
 
 \bi \item[] \it Throughout this section, we  assume  $A = A_\pi$ is a closed translation-invariant (necessarily unital) subalgebra of $B(G)$ that contains $A_\F$.   \rm
 
 \ei  We shall begin with an examination of the structure of $\Delta(A)$, a compact semitopological $*$-semigroup by Corollary \ref{SpectrumSemigpCor}, and will then employ our analysis in a proof of the converse to Proposition   \ref{PWAffineCBProp} in the case that $\Delta(A) = G \cupdot G^{ap}$.  For example, we  characterize all completely bounded  homomorphisms $\vp:B(G) \ra B(H)$ when $G$ is a Euclidean- or $p$-adic-motion group.  In the process, we aim to illustrate some of the results from the previous sections prior to turning our focus to the study of completely positive/completely contractive homomorphisms in  Section 4.  
 
 As shown by Runde and Spronk \cite[Proposition 2.1]{Run-Spr0},  the algebra $A_{\cal F}$ is the closed linear span in $B(G)$ of the coefficient functions of the family $\cal F$ of finite-dimensional unitary representations of $G$, so  $A_{\cal F} = A_{\pi_{\cal F}}$ for a representation $\pi_{\cal F}$ of $G$ (specifically, $\pi_{\cal F} $ is the direct sum of all finite dimensional irreducible representations of $G$) and $j_{\delta^{ap}}: v \mapsto v \circ \delta^{ap}$ is a completely isometric algebra isomorphism of $A(G^{ap})$ onto  $A_{\cal F} = A(G^{ap}) \circ \delta^{ap}$; also see  \cite[(2.27)]{Eym} and \cite[Corollary  2.2.5]{Kan-Lau}. Letting $A_{\cal PIF}$ denote the unique vector space complement of $A_\F$ in $B(G)$---see \cite[Theorem 3.18]{Ars} and \cite[Remark 2.5]{Run-Spr}---$A_{\cal PIF}$ is a closed translation-invariant ideal in $B(G)$ by \cite[Theorem 2.3]{Run-Spr0}. Since we are assuming that $A_\F$ is contained in $A$, it follows that $A_{\F^\perp}: = A_{\cal PIF} \cap A$ is a closed translation-invariant ideal in $A$ and we have the $\ell^1$-direct sum decomposition of $A$,  $A = A_{\F^\perp} \oplus_1 A_\F$ \cite[Theorem 3.18 and Corollary 3.13]{Ars}. Hence, each $w$ in $A$ takes the form $w = u + v \circ \delta^{ap}$ for some unique $u \in A_{\F^\perp}$ and $v \in A(G^{ap})$.  If $A = A_\F$, the results below are trivial,  so \it in this section   we will always assume  that $G$ is non-compact and   $A_{\F^\perp}$ is nontrivial. \rm
 
  Let $e_{ap}$ denote the identity in $G^{ap}$ and let $\Delta(A)_+ = \{ s \in \Delta(A): s \geq 0\}$, viewing  $\Delta(A)$ as a subset of  $ VN_\pi$.
  Define  $$\psi: G^{ap} \ra \Delta(A) \bs \pi(G): a \mapsto \psi_a \ \ {\rm by \ putting} \ \ \psi_a(u + v\circ \delta^{ap}) = v(a).$$  
(Observe that $\psi$ maps into $\Delta(A)$ since $A_{\F^\perp}$ is an ideal in $A$;  $A_{\F^\perp}$ is nontrivial and translation-invariant so it vanishes at no point of $G$,  from which it follows that  $\psi$ maps into $\Delta(A) \bs \pi(G)$.)

   \bt  \label{Delta(A) and ap-cpctn Thm} The following statements hold:
   \bi \item[(i)] The compactification homomorphism $\delta^{ap}: G \ra G^{ap}$ determines a continuous $*$-epimorphism $$\widetilde{\delta^{ap}}: \Delta(A) \ra G^{ap} \ \  such \ that \ \ \tdap(\pi(g)) =  \dap(g) \ \ and \ \ \tdap (s) = \eap$$ for any $g \in G$ and $s \in \Delta(A)_+$. 
    \item[(ii)] The map $\psi: \gap \ra \Delta(A) \bs \pi(G)$ is a topological isomorphism of $\gap$ onto its image $\psi(\gap)$, which is a compact ideal in $\Delta(A)$. Moreover,  for any $a \in \gap$ and $x \in \Delta(A)$, 
   \beq  \label{Delta(A) and ap-cpctn Thm Eqn} \psi_a x = \psi_{a \tdap(x)} \quad   and  \quad  x\psi_a = \psi_{\tdap(x)a}; \eeq  more concisely,  identifying $\gap$ and $\psi(\gap)$ we have     \beq  \label{Delta(A) and ap-cpctn Thm Eqn2} a x = a \tdap(x) \qquad and \qquad xa= \tdap(x)a.  \eeq 
   \ei 
   
   \et 
   
   \begin{proof} (i) By assumption $A$ contains $A_\F = A(\gap) \circ \dap$, so $j_{\dap}$ maps $A(\gap)$ into $A = A_\pi$ and has dual map $\gamma_{ap}: = j_{\dap}^*: VN_\pi \ra VN(\gap)$ satisfying $\gamma_{ap}(\pi(g)) = \lambda_{\gap}(\dap(g)) = \dap(g)$. As noted in the proof of Proposition \ref{CPHomProp}, it follows that $\gamma_{ap}$ is a continuous $*$-homomorphism. Moreover, since $j_{\dap}:A(\gap) \ra A$ is a homomorphism mapping $1_{\gap}$ to $1_G$, $\gamma_{ap}$ maps $\Delta(A)$ into $\Delta(A(\gap)) =\gap$. Hence, 
   $$ \tdap: = \gamma_{ap}{\large |}_{\Delta(A)}: \Delta(A) \ra  \gap$$
 is a continuous closed-range $*$-homomorphism such that $\tdap(\pi(g))= \dap(g)$ ($g \in G$); hence, $\tdap$ is surjective. For $x \in \Delta(A)$ and $v \in A(\gap)$, observe that 
 \beq  \label{Delta(A) and ap-cpctn Thm Eqn3}  v(\tdap(x)) = \l \tdap(x), v \r = \l x, v \circ \dap \r. \eeq     
 Let $s\in \Delta(A)_+$. Then, since $\gamma_{ap}$ is a $*$-homomorphism, $\tdap(s) \geq 0$ and $\tdap(s) = \tdap(s^*) = \tdap(s)^*$, i.e., $\tdap(s) = \tdap(s)^{-1}$ in $\gap$. Thus, $\tdap(s)^2 = \eap = \eap^2$, meaning that $\tdap(s) = \eap$, the unique positive square root of $\eap$ in $VN(\gap)$.   
 
 (ii)   The first part of the proof is based on the proof of \cite[Prop. 2.9.2]{Kan-Lau}: From the Gelfand theory, $Q^*: \vp \mapsto \vp \circ Q$ is a homeomorphism of $\Delta(A(G^{ap}) \circ \delta^{ap}) = \Delta (A / A_{\F^\perp})$ and $h(A_{\F^\perp}) \subseteq \Delta(A)$, the hull of the ideal $A_{\F^\perp}$,  where $Q: A \ra A(G^{ap}) \circ \delta^{ap}$ is the quotient map. Since $\kappa_{ap} = (j_{\delta^{ap}})^{-1}: A(G^{ap}) \circ \delta^{ap} \ra A(G^{ap})$ is an isometric algebra isomorphism, $\kappa_{ap}^*$ maps $G^{ap} = \Delta(A(G^{ap}))$ homeomorphically onto $\Delta(A(G^{ap}) \circ \delta^{ap})$. Hence, $Q^* \circ \kappa_{ap}^*$ is a homeomorphism of $G^{ap}$ onto $h(A_{\F^\perp}) \subseteq \ \Delta(A)$, and it is easy to see that $Q^* \circ k_{ap}^* = \psi$. 
 
 To see that $\psi$ is a homomorphism and (\ref{Delta(A) and ap-cpctn Thm Eqn}) holds, we  employ the Arens product formulation of multiplication  in $\Delta(A)$  from Proposition \ref{ArensProductProdProp}: Let $a, b \in G^{ap}$, $w = u +  v\circ \delta^{ap}  \in A$.  For $s \in G$,   
\beqs  \psi_b \c ( u+ v\circ \delta^{ap}) (s) & = & \psi_b( (u+ v\circ \delta^{ap})\c s) = \psi_b( u \c s +  (v \c \delta^{ap}(s)) \circ \delta^{ap}) \\
& = &  (v \c \delta^{ap}(s)) (b) = v(\delta^{ap}(s)b) = (b \c v) \circ \delta^{ap}(s).     \eeqs
Since $b \c v \in A(G^{ap})$, we obtain
\beqs \psi_a \psi_b (u+ v \circ \delta^{ap}) & = &   \psi_a( \psi_b \c ( u+ v \circ \delta^{ap})) = \psi_a ((b \c v) \circ \delta^{ap})  \\ 
& = & (b\c v)(a) = v(ab) = \psi_{ab}(u+ v \circ \delta^{ap}),\eeqs
as needed. For $x \in \Delta(A)$, 
\beqs \l x \psi_a, u+  v \circ \delta^{ap}\r & = & \l x, \psi_a \c (u+ v \circ \delta^{ap})\r= \l x,   (a \c v)\circ \delta^{ap}\r \\ 
& = &  (a \c v)(\tdap(x))  = v(\tdap(x)a) = \l  \psi_{\tdap(x) a} , u+  v \circ \delta^{ap}\r,  \eeqs
using (\ref{Delta(A) and ap-cpctn Thm Eqn3}). Hence, $x \psi_a = \psi_{\tdap(x)a} $ and a symmetric argument shows  $\psi_a \, x = \psi_{a \tdap(x)}$.       
   \end{proof} 
   
The following corollary  includes simple new proofs of some of the main results from \cite{Wal2}, (ibid. Theorem 2 and Proposition 8). More generally, we prove the result for any Banach subalgebra $A = A_\pi$ of $B(G)$ that contains $A_{\F}$. \rm  
 
 \bc \label{Walter min idempotent etc Corollary}  The following statements hold: 
 \bi \item[(i)]  $z_F: = \psi_{e_{ap}}$ is a central idempotent in  $\Delta(A)$  such that for any $s \in \Delta(A)_+$, $s z_F = z_F$; so, $z_F = \min(\Delta(A)_+)$. 
 \item[(ii)]  Define $\theta_F: G \ra \bD_{z_F}$ by $\theta_F(g) = \pi(g) z_F$. Then $\psi$  is a compactification isomorphism of $(\gap, \dap)$ onto $(\bD_{z_F}, \theta_F)$; so,  $(\bD_{z_F}, \theta_F)$  is the almost periodic compactification of $G$ and $\bD_{z_F} = \psi(G^{ap})$ is a compact ideal in $\Delta(A)$.
\item[(iii)] $A_{\F} = z_F \c A$ and $A_{\F^\perp} = (\pi(e_G) - z_F)  \c A$.
\item[(iv)]  $G$ is compact if and only if $e_G$ is the only (central) idempotent in $\Delta(B(G))$.  
 \ei 
 \ec 
 
 \begin{proof}  (i) Let $x \in \Delta(A)$, $s \in\Delta(A)_+$. Identifying $\gap$ and $\psi(\gap)$, (\ref{Delta(A) and ap-cpctn Thm Eqn2}) and Theorem \ref{Delta(A) and ap-cpctn Thm}(i) give $$\eap x = \eap \tdap(x) = \tdap(x) = \tdap(x) \eap = x \eap \quad {\rm and} \quad s \eap = \tdap(s) \eap = \eap.$$ We can thus take $z_F= \eap (=\psi_{\eap})$. 
 
 (ii)  By Proposition \ref{SpectralSubgroupsProp}, $\psi(\gap) \subseteq \bD_{\psi(\eap)} = \bD_{z_F}$ and for any $x \in \bD_{z_F}$,   $$x = xz_F = x \psi_{\eap} = \psi_{\tdap(x) \eap} = \psi_{\tdap(x)} \in \psi(\gap).$$
Hence, $\psi:\gap\ra \psi(\gap) = \bD_{z_F}$ is a topological isomorphism and for any $g \in G$, 
$$\theta_F(g) = \pi(g) z_F = \pi(g) \psi_{\eap} = \psi_{\tdap(\pi(g))\eap} = \psi_{\dap(g)} = \psi \circ \dap(g).$$ 
It follows that $\theta_F$ is a continuous dense-range homomorphism and $\psi$ is an isomorphism of the compactifications $(\gap, \dap)$ and $(\bD_{z_F}, \theta_F)$. 

(iii)  Let $w = u + v \circ \dap \in A$ where $u \in A_{\F^\perp}$ and $v \in A(\gap)$.  For $g \in G$, 
\beqs  z_F \c w(g) = \l \pi(g), \psi_{\eap}  \c w \r = \l \pi(g) \psi_{\eap} , w \r = \l \psi_{\dap(g)}, u + v \circ \dap \r =  v \circ \dap(g).      \eeqs
Hence, $z_F \c A = A(\gap) \circ \dap = A_\F$ and $(\pi(e_G)- z_F) \c A = A_{\F^\perp}$. 

(iv) If $G$ is non-compact, we have noted that $\psi$ maps into $\DBG \bs \omega_G(G)$, so   $z_F = \psi_{\eap}$ is a central idempotent in $\DBG$ that is different from $e_G$.  If $G$ is compact, then $\DBG = \Delta(A(G)) = G$, so the converse implication is trivial. 
 \end{proof}

 \br \rm From the theory of compact right topological semigroups, $S= \Delta(A)$ has a smallest two-sided ideal $K(S)$ \cite[Theorem II.2.2]{Ber-Jun-Mil1}.  By Corollary \ref{Walter min idempotent etc Corollary}, $\bD_{z_F} = \psi(\gap)$ is a two-sided  ideal in $S$. Moreover, if $L$ is any left ideal in $S$ and $x \in L$,  $\psi(\gap) =\{\psi_{a\tdap(x)}: a \in \gap\} = \psi(\gap)x$, which is contained in $L$. Hence, Theorem \ref{Delta(A) and ap-cpctn Thm} allows us to identify $K(S)$ as $\psi(\gap) = \bD_{z_F}$, the almost periodic compactification of $G$. 
 \er

As above, express $A = A_\pi$ as $A = A_{\F^\perp} \oplus_1 A_\F$. Then  $\pi (G) \subseteq \Delta(A)$ and taking $\pi_{\F^\perp}$ such that $A_{\pi_{\F^\perp}} = A_{\F^\perp}$, $\pi_{\F^\perp}(G) \subseteq \Delta(A_{\F^\perp})$. 

\bc \label{Delta(A) = G cupdot Gap Cor}   Suppose that  $\Delta (A_{\F^\perp}) = \pi_{\F^\perp}(G)$ and $A$ is a regular algebra of functions on $G$. Then 
$$\Delta(A) = \pi(G) \cupdot \psi(G^{ap})$$ 
where $\pi: G \ra \pi(G) $ and $\psi: \gap \ra  \psi(\gap)$ are topological  isomorphisms.  Moreover,  $\psi(\gap)= \bD_{z_F}$ is a compact ideal in $\Delta(A)$ satisfying 
$$\psi_a \pi(g) = \psi_{a \dap(g)} \quad and \quad    \pi(g) \psi_a = \psi_{\dap(g)a}    \qquad (g \in G, \ a \in \gap).$$ 
 More concisely, if we identify $G$ with $\pi(G)$ and $\gap$ with $\psi(\gap)$, we have  $$\Delta(A) = G \cupdot \gap$$
 where $G$ and $\gap$ are open and compact subgroups of $\Delta(A)$ respectively, and $\gap$ is an ideal in $\Delta(A)$ satisfying 
 $$ag = a \dap(g) \quad { and} \quad   ga = \dap(g) a   \qquad (g \in G, \ a \in \gap).$$
  \ec 
  
  \begin{proof}   We are assuming that $A= A_\pi$ is a regular algebra of functions on $G$, so $\pi:G \ra \pi(G)$ is a topological isomorphism. Let $\vp \in \Delta(A) \bs \psi(\gap)$. From the proof of Theorem \ref{Delta(A) and ap-cpctn Thm}(ii), $\psi(\gap) = h(A_{\F^\perp}) \subseteq \Delta(A) \bs \pi(G)$, the hull of the ideal $A_{\F^\perp}$, so $\vp\large{|}_{A_{\F^\perp}} \in \Delta (A_{\F^\perp})$.   Let $g \in G$ be such that  $\vp\large{|}_{A_{\F^\perp}}  = 
\pi_{\F^\perp}(g)$, and choose $u \in A_{\F^\perp}$ such that $\vp(u) = u(g) =1$. For any $w \in A$, $uw \in A_{\F^\perp}$, so  $\vp(w) = \vp(u) \vp(w) = \vp (uw) = u(g) w(g) = \pi(g)(w)$; hence, $\vp = \pi(g)$. Thus, $\Delta(A) = \pi(G) \cupdot \psi(\gap)$ and the remaining statements  are immediate from Theorem \ref{Delta(A) and ap-cpctn Thm}.  
  \end{proof}  
  
\bex   \label{AFG Type Examples}  \rm  (a)  \it Introduction of $\AFG$, an accessible test subalgebra of $B(G)$: \rm    Let $G$ be any non-compact locally compact group and let 
$$\AFG:= A(G) \oplus_1 A_{\F} =A(G) \oplus_1 A(\gap) \circ \dap.$$
Since $A(G)$ is an ideal in $B(G)$, $A = \AFG$ is a closed translation-invariant unital subalgebra of $B(G)$. By \cite[Theorem 3.18]{Ars},  we must have $A_{\F^\perp} = A(G) = A_{\lambda_G}$, where $\lambda_G = \pi_{\F^\perp}$ is the left regular representation of $G$;  letting $\pif := \lambda_G \oplus \pi_{\cal F}$,   $\AFG = A_\pif$. Since $\Delta(A(G)) = \lambda_G(G)$, Corollary \ref{Delta(A) = G cupdot Gap Cor} applies and we obtain  $$\Delta(\AFG) = \pif(G) \cupdot \psi(\gap) = G \cupdot \gap$$ 
with the product  and identifications as defined in Corollary \ref{Delta(A) = G cupdot Gap Cor}.


With Corollary \ref{Full characterization of cc and cp in special cases Cor}---specifically, see  Remark \ref{CP and CC Homom G cup Gap Remarks}---we shall see that while $\AFG$ and its spectrum are tractable, they are  nonetheless intricate enough to illustrate how completely positive/completely contractive homomorphisms of (subalgebras of) $B(G)$ can be significantly more complicated than the same types of  homomorphisms on  $A(G)$.   When $G$ is abelian with dual group $\widehat{G}$,  $L^1(G)^\wedge$ is  $A(\widehat{G})$ and, as shown by Runde and Spronk \cite[p. 679]{Run-Spr0}, $\ell^1(G)^\wedge$ is the subalgebra $B(\widehat{G}) \cap AP(\widehat{G}) = A_{\cal F}(\widehat{G})$, where $\mu \mapsto \hat{\mu}$ is the Fourier--Stieltjes transform on $M(G)$. Thus,  
  $\AFG $ is dual to the subalgebra $M_a(G) \oplus_{\ell^1} M_d(G)$ of $M(G)$, where $M_a(G) = L^1(G)$ is the space of measures that are absolutely continuous with respect to Haar measure, and $M_d(G) = \ell^1(G)$ is the space of discrete measures.  Letting $M_{cs}(G)$ denote the continuous measures that are singular with respect to Haar measure,  the most tractable part of $M(G) =M_a(G)  \oplus_{\ell^1} M_{cs}(G) \oplus_{\ell^1} M_d(G)$ (e.g., see \cite[Ch. 3]{Dal})  is certainly the  Banach algebra of measures whose continuous parts are absolutely continuous, $M_{\mathbb F}(G):= M_a(G) \oplus_{\ell^1} M_d(G) = L^1(G) \oplus_{\ell^1} \ell^1(G) $.  The (completely) positive and contractive homomorphisms  $\vp : L^1(G) \ra M(H)$ have been characterized for arbitrary locally compact groups \cite{Gre, Sto1}, and studying the positive and contractive homomorphisms $\vp: M_{\mathbb F}(G) \ra M(H)$  will be a good entry point to an investigation of dual versions of some of the problems discussed in this paper. 
  
  \smallskip 
  
  \noindent (b) \it $B(G)$ for  a $p$-adic-motion group $G$: \rm   More generally, we  consider the locally compact groups studied by Runde and Spronk in \cite{Run-Spr}.   Let $G= N \rtimes K$ be a semidirect product of a compact group $K$ acting on a non-compact locally compact abelian group $N$ such that properties (1) and (2) on page 196 of \cite{Run-Spr} (also see Section 2.9 of \cite{Kan-Lau}) are satisified.  Then,  \beq \label{FS Decomp for Fell groups} B(G) = A(G)  \oplus_1 A(K) \circ q \eeq
\cite[Prop. 2.1]{Run-Spr} or \cite[Prop. 2.9.1]{Kan-Lau}.  Letting $q: G = N \rtimes K \ra K$ be the quotient map, we observe that 
\bi \item[(i)] $A_{\cal F} = A(K) \circ q$; and 
\item[(ii)] $(K,q) = (\gap, \dap)$,  the almost periodic compactification of $G$.  
\ei 
(These statements are implicit elsewhere, but for the convenience of the reader we provide some explanation:  Note that $j_q$ is continuous (isometric, in fact) with respect to uniform norms,  every continuous function on $K$ is almost periodic since $K$ is compact, and $\ell_x(j_q(v)) = j_q(\ell_{q(x)} v)$ for $v \in C(K)$; hence $j_q$ maps $C(K)$ into $AP(G)$ and it follows that $A(K)\circ q = j_q(A(K))$ is contained in $B(G) \cap AP(G) = A_{\cal F}$. On the other hand, $w \in A_{\cal F}$ can be written as $w= u + v \circ q$ for some $u \in A(G)$ and $v \in A(K)$ by (\ref{FS Decomp for Fell groups}). But $v \circ q \in A_{\cal F}$, so $u = w - v\circ q \in A(G) \cap A_{\cal F} = \{0\}$, whence $w = v \circ q$. This establishes (i).  Thus,  $j_q$ and $j_{\delta^{ap}}$ are isometric algebra isomorphisms mapping $A(K)$ and $A(G^{ap})$, respectfully, onto $A_{\cal F}$.   Therefore, $\vp:= (j_q)^{-1} \circ  j_{\delta^{ap}}$ is an isometric algebra isomorphism of $A(G^{ap})$ onto $A(K)$,  $\theta= \vp^*{\large |}_K $ maps $K = \Delta(A(K))$ homeomorphically onto $\Delta(A(G^{ap}))= G^{ap}$, and it is easy to see that $\theta \circ q = \delta^{ap}$. It  automatically follows that $\theta$ is a homomorphism  \cite{Ber-Jun-Mil} and (ii) is  established.)

We conclude that  $B(G) = A(G)  \oplus_1 A(K) \circ q= \AFG$ in this case and therefore, as above, $\DBG =  G \cupdot K$ with $K = \gap$ a compact  ideal in $\DBG$ satisfying  
$$ a x = a q(x) \ \ {\rm and} \ \ x a = q(x) a \qquad 
 (a \in K, \ x \in G).$$   
 
 For a concrete example, let $p$ be a prime number, $\mathbb{Q}_p$ the locally compact topological field of $p$-adic numbers with valuation $|\c|_p$, $\mathbb{O}_p = \{ r \in \mathbb{Q}_p: |r|_p \leq 1\}$ the compact open subring of $p$-adic integers. (Section III.6 of \cite{FelDor} contains a concise introduction to these objects.)   Letting $K = SL(n, \O_p)$, the compact multiplicative group of $n\times n$ matrices $A$ over $\O_p$ with $|\det(A)|_p=1$, act on the vector space $N = \Q_p^n$ by matrix multiplication, $G_{p,n}= N\rtimes K$---called the ``$n$th rigid $p$-adic motion group" in \cite{Run-Spr}---is an example of a locally compact group satisfying these two properties \cite[Example 2.9.3]{Kan-Lau}.  When $n=1$, $G_{p,n}$ is the Fell group $\Q_p \rtimes \mathbb{T}_p$, where $\mathbb{T}_p= \{ r \in \Q_p: |r|_p=1\}$.  
 
 \smallskip 
 
 \noindent (c) \it $B(G)$ for a Euclidean motion group $G$:    \rm 
 Let $G = {\mathbb R}^n \rtimes SO(n)$ with $n\geq 2$, where $SO(n)$ acts on ${\mathbb R}^n$ by matrix multiplication (rotation). Letting $q: G \ra SO(n)$ be the quotient map and $B_0(G) = B(G) \cap C_0(G)$ the Rajchman algebra of $G$,  $\Delta(B_0(G)) = G$ and $B(G) = B_0(G) \oplus_1 B(SO(n)) \circ q$ by \cite[Theorem 4.1]{Kan-Lau-Ulg}.   Since $G$ is non-compact, $C_0(G) \cap AP(G) = \{0\}$ and the argument given above in (b) shows that $A_{\F} = A(SO(n)) \circ q$ and $(SO(n), q ) = (\gap, \dap)$. Hence, $A_{\F^\perp} = A_{\cal PIF}$, the unique vector space complement of $A_\F$ in $B(G)$, is $B_0(G)$ and $\Delta(A_{\F^\perp}) = G (= \pi_{\F^\perp}(G)$).   By Corollary \ref{Delta(A) = G cupdot Gap Cor}, $\Delta(B(G)) = G \cupdot SO(n)$ (also shown in \cite[Theorem 4.1(iii)]{Kan-Lau-Ulg}) and, moreover, $SO(n)$ is a compact ideal in $\Delta(B(G))$ satisfying 
 $$ a x = a q(x) \ \ {\rm and} \ \ x a = q(x) a \qquad 
 (a \in SO(n), \ x \in G).$$

  \smallskip 
 
 \noindent (d) \it Spine examples:  \rm  
 In \cite{Il-Sp2}, the authors introduced the spine $A^*(G)$ of $B(G)$ and studied  $G^*= \Delta(A^*(G))$, the spine compactification of $G$.      They showed that  
 $$A^*(G) = A(G)  \oplus_1 A_{\tau_{ap}}(G) = A(G) \oplus_1 A(\gap) \circ \dap = \AFG$$  when  $G$ contains a compact normal subgroup $K$ for which $G/K$ is topologically isomorphic to  $\mathbb{R}$ or $\mathbb{Z}$  \cite[Proposition 6.2]{Il-Sp2};   $G$ is the $p$-adic field $\mathbb{Q}_p$ for some prime number $p$ \cite[Proposition 6.7]{Il-Sp2}; or $G$  is  any minimally weakly almost periodic group such as $G= SL_2(\mathbb{R})$ or $G=\mathbb{R}^n \rtimes SO(n)$  by Theorem 6.5 and page 297 paragraph 1 of \cite{Il-Sp2}.  Thus, in all of these cases $\Delta(A^*(G)) = G^*  = G \cupdot G^{ap}$ with multiplication as described as in Corollary  \ref{Delta(A) = G cupdot Gap Cor}.  (This description of $G^* = \Delta(A^*(G))$ is explicitly derived in \cite{Il-Sp2}.)   
\eex

 \bp \label{Singular Hom on AFG Prop} Let  $A$ be a regular algebra of functions on $G$ such that $\Delta(A_{\F^\perp}) = \pi_{\F^\perp}(G)$ and  $\vp: A \ra B(H)$.  Then $\vp$ is  a  singular completely bounded/ completely contractive/ completely positive homomorphism if and only if $\vp = j_\alpha$  where  $\alpha: Y  \ra  \bD_{z_F}$  is a continuous piecewise-affine/ affine/ homomorphic map. 
\ep 

\begin{proof}  By Corollary \ref{Delta(A) = G cupdot Gap Cor}, $\Delta(A) = \pi(G) \cupdot \psi(\gap) = G \cupdot \bD_{z_F}$.    Suppose that $\vp$ is a singular homomorphism, so $\vp = \ja$ where   $\alpha: Y \ra   \psi(\gap)$  is  continuous on $Y \in \Omega(H)$ (Proposition \ref{FSAlgHomGeneralProp}).  The map 
$\vp_{ap} := \vp \circ j_{\delta^{ap}}: A(G^{ap}) \ra B(H)$ is a completely bounded/completely contractive/completely positive homomorphism and  $\vp_{ap}= j_{\alpha_{ap}}$ where $\alpha_{ap} = \psi^{-1} \circ \alpha:Y \ra G^{ap}$. To see this, let $v \in A(G^{ap})$. For $y \in Y$, $\alpha(y) = \psi_a$ for some $a \in G^{ap}$, and therefore $\alpha_{ap}(y) = a$. Hence, 
$$\vp_{ap}(v) (y) = j_\alpha(v\circ {\delta^{ap}})(y) = \l \alpha(y), v \circ \delta^{ap} \r = \l \psi_a, v\circ \delta^{ap}\r = v(a)  = j_{\alpha_{ap}} v(y).$$ 
For $y \in H \bs Y$, $\vp_{ap}(v)(y) = j_\alpha(v\circ \delta^{ap})(y) = 0 = j_{\alpha_{ap}}v(y)$, as needed.  Since $G^{ap}$ is amenable, by \cite[Theorem 3.7]{Il-Sp} (also Proposition \ref{PhamProp5.8} and Remark \ref{PhamProp 5.8 Remark}.2), $\alpha_{ap}$ is piecewise-affine/affine/homomorphic; as $\psi$ is a topological isomorphism of $G^{ap}$ onto $\psi(\gap) = \bD_{z_F}$, $\alpha = \psi \circ \alpha_{ap}$ is the same.  The converse statement is an immediate consequence of Proposition \ref{PWAffineCBProp}/Proposition \ref{CCProp}/Corollary \ref{CPHomCor}.
 \end{proof}  

 \bt \label{AFG Main Homom Thm} Suppose that  $G$ is   amenable,  $A$ contains $A(G)$ and $\Delta(A_{\F^\perp}) = \pi_{\F^\perp}(G)$.  Then $\vp:A \ra B(H)$ is a completely bounded homomorphism if and only if $\vp= \ja$ for a continuous pw$^2$-affine map $\alpha:Y \ra \Delta(A)$.   \et

\begin{proof} By Corollary \ref{Lebesgue Decomp Cor}(i), $\vp$ has a unique Lebesgue decomposition $\vp = \vp_a + \vp_s = j_{\alpha_a} + j_{\alpha_s}$ such that $\alpha_a: Y_a \ra G$ is continuous and piecewise affine.  Since $\vp_s = \vp - \vp_a$ is completely bounded and singular, $\alpha_s: Y_s \ra \bD_{z_F}$ is continuous and piecewise affine by  Proposition \ref{Singular Hom on AFG Prop}. Therefore, $\alpha= \alpha_a \cupdot \alpha_s$ is a continuous pw$^2$-affine map. The converse is a special case of Proposition \ref{PWAffineCBProp}.
\end{proof}

\br \rm  1. Proposition \ref{Singular Hom on AFG Prop} and Theorem \ref{AFG Main Homom Thm} apply to all of the examples of $A$ considered in Example \ref{AFG Type Examples}. For instance, when $G$ is a Euclidean- or $p$-adic-motion group, $\vp=\ja:B(G) \ra B(H)$ is a completely  bounded homomorphism if and only if $\alpha$ is a continuous pw$^2$-affine map. However,  the methods used to prove Theorem \ref{AFG Main Homom Thm}---partly based on the main theorem from \cite{Il-Sp}---cannot be used to characterize the completely positive and completely contractive homomorphisms $\vp: A \ra B(H)$ in these cases; these homomorphisms will be fully described in Section 4.   \\
2. With Example \ref{Rudin Ex as pw2 affine}, we showed that it is easy to obtain examples of completely bounded homomorphisms $\vp: A \ra B(H)$  such that  $\alpha$ does not map into either $G$ or $\bD_{z_F}$.   In Section 4 we will give completely positive and completely contractive examples for which this is the case.

\er

\section{Completely positive and completely contractive homomorphisms: necessary conditions}

 \bi \item[] \it Unless stated otherwise, throughout this section $A=A_\pi$ is a  closed unital translation-invariant subalgebra of $B(G)$.  \rm 
  \ei

 \blem \label{PositiveHomLem}  Let $\vp:A \ra B(H)$ be a positive homomorphism. Then $\vp$ is automatically contractive and $\vp = \ja$ for some continuous map $\alpha:H_\alpha \ra \Delta(A)$ where $H_\alpha$ is an open subgroup of $H$ and $\alpha(h^{-1}) = \alpha(h)^*$ for each $h \in H_\alpha$.  Moreover, if $\alpha(e_H) \in \bD_p$ for some projection $p \in \Delta(A)$, then $\alpha(e_H) = p$. 
 \elem

 \begin{proof} By Proposition \ref{FSAlgHomGeneralProp}, $\vp = \ja$ where $Y = \{ h \in H: \vp^*(h) \neq 0\}$ and $\alpha:Y \ra \Delta(A)$ is the restriction of $\vp^*$ to $Y$. Since $\vp$ is positive, $1_Y = \vp(1_G)$ is positive definite, whence $Y=H_\alpha$ is an open subgroup of $H$ by Corollary \ref{Il-SpThm2.1Cor}, and by  Theorem \ref{CPBasicThm}(a),  $\| \vp\| = \|\vp^*\| =\|\vp^*(e_H)\|= 1$ since $\vp^*(e_H) \in \Delta(A)$. For each $x \in W^*(H)$, $\vp^*(x^*) = \vp^*(x)^*$---also by Theorem \ref{CPBasicThm}(a)---so for each $h \in H_\alpha$, $\alpha(h^{-1}) = \alpha(h)^*$. Finally, suppose that $\alpha(e_H) \in \bD_p$. Since $e_H \geq 0$, $\alpha(e_H) = \vp^*(e_H) \geq 0$ and, under our hypothesis, $\alpha(e_H)^2= \alpha(e_H) \alpha(e_H)^* = p$, so $\alpha(e_H) = p^{1/2}= p$. 
 \end{proof}  


In the following theorem and its corollaries, we assume $\vp:A \ra B(H)$ is a positive homomorphism, which---by Lemma \ref{PositiveHomLem}---takes the form $\vp = \ja$ where $\alpha: H_\alpha \ra \Delta(A)$ is a continuous map on an open subgroup $H_\alpha$ of $H$ satisfying $\alpha(h^{-1}) = \alpha(h)^*$ for each $h \in H_\alpha$.  
 
\bt \label{CPMainThm}  Suppose that $\alpha(e_H)$ belongs to a subgroup of $\Delta(A)$, say $\alpha(e_H) \in \bD_p$, and  \bi \item[(i)]  $\vp = \ja$ is completely positive, or 
\item[(ii)] $\vp = \ja$ is 4-positive and $\alpha(H_\alpha)$ is contained in $p \Delta(A)p$. 
\ei Then $\alpha(e_H) = p$, 
$$H_p = \{ h \in H_\alpha: \alpha(h) \in \bD_p\}$$ 
is a subgroup of $H_\alpha$,  $\alpha:H_p \ra \bD_p$ is a continuous homomorphism, and $$\alpha (hk) = \alpha(h) \alpha(k), \quad \alpha(kh) = \alpha(k) \alpha(h) \ \ and  \ \ \alpha(h^{-1}) = \alpha(h)^* $$  
 for all $h \in H_\alpha$ and $k \in H_p$. When $\vp$ is completely positive, $\alpha(H_\alpha)$ is contained in $p \Delta(A)p$.
\et 


\begin{proof}  Suppose that $\vp$ is completely positive. By Lemma \ref{PositiveHomLem},  for $k \in H_p$,  $\alpha(k^{-1}) = \alpha(k)^* \in \bD_p$, so $k^{-1} \in H_p$; moreover $\alpha(e_H) = p$, so  $\vp^*(k)^* \vp^*(k) = p = \vp^*(k^*k)$.  
Since  $\vp$ is contractive,     
$$  \alpha(hk) = \vp^*(hk) = \vp^*(h) \vp^*(k) = \alpha(h) \alpha(k)       $$
whenever $h \in H_\alpha$ and $k \in H_p$ (and, similarly, $\alpha(kh) = \alpha(k) \alpha(h)$), by Theorem \ref{CPBasicThm}(c). In particular, $hk \in H_p$ whenever $h, k \in H_p$, so $H_p$ is a subgroup of $H_\alpha$ and $\alpha:H_p \ra \bD_p$ is a homomorphism.  Observe that for  any $h \in H_\alpha$, $\alpha(h) = \alpha(he_H) = \alpha(h) \alpha(e_H) = \alpha(h)p$ and $\alpha(h) = p\alpha(h)$, so $\alpha$ maps $H_\alpha$ into $p \Delta(A) p$ in this case. If condition (ii) holds, then $\vp^*$ maps  $H$ into the von Neumann subalgebra   $p VN_\pi p$ of $VN_\pi$. By weak$^*$-continuity of $\vp^*$ and weak$^*$-density  of the linear span of $H$ in $W^*(H)$, $\vp^*$ maps $W^*(H)$ into $p VN_\pi p$ and $\vp^*(e_H) = p$, the identity in $p VN_\pi p$. Hence, Theorem  \ref{CPBasicThm}(c) can also be applied in this case, and the conclusion follows as before.   
\end{proof} 

As an immediate consequence of Proposition \ref{SpectralSubgroupsProp}, Corollary \ref{CPHomCor} and Theorem \ref{CPMainThm},   we have the following corollary via the chain of implications $(i) \Rightarrow (iii) \Rightarrow (ii) \Rightarrow (i)$.  

\bc \label{CPCor1} Let  $\vp = \ja:A \ra B(H)$ be a homomorphism. The following statements are equivalent: \bi \item[(i)]   $\vp$ is 4-positive and $\alpha(Y)$ is contained in a subgroup of $\Delta(A)$;
\item[(ii)]   $\vp$ is completely positive and $\alpha(Y)$ is contained in a subgroup of $\Delta(A)$; 
\item[(iii)] $Y$ is an open subgroup $H_\alpha$ of $H$ and $\alpha:H_\alpha \ra \Delta(A)$ is a continuous homomorphism (into $\bD_p$ for some projection $p$ in $\Delta(A)$).
\ei  
 \ec


 
 The following result improves Proposition \ref{Singular Hom on AFG Prop} in the completely positive case. 
 
 \bc \label{CP Cor z_F} Suppose that $A= A_\pi$ contains $A_\F$ and let $\vp = \ja: A \ra B(H)$ be a homomorphism. Then the following statements are equivalent: 
 \bi \item[(i)] $\vp$ is completely positive and $\alpha(e_H) \in \bD_{z_F}$; 
 \item[(ii)] $Y$ is an open subgroup of $H$ and $\alpha: Y \ra \bD_{z_F}$ is a continuous homomorphism;
 \item[(iii)] $\vp$ is completely positive and $A_{\F^\perp}$ is contained in the kernel of  $\vp$. 
 \ei 
 \ec

 \begin{proof}  By Corollary \ref{Walter min idempotent etc Corollary}, $A_\F = z_F \c A$ and $A_{\F^\perp} = (\pi(e_G) - z_F) \c A$. \\
 (i) implies (ii): By Lemma \ref{PositiveHomLem} and Theorem \ref{CPMainThm}, $Y=H_\alpha$ is an open subgroup of $H$, $\alpha(e_H) = z_F$ and $\alpha(H_\alpha)$ is contained in $z_F \Delta(A) z_F$. By  Corollary \ref{Walter min idempotent etc Corollary},  $\bD_{z_F}$ is an ideal in $\Delta(A)$, so $\alpha(H_\alpha) \subseteq \bD_{z_F}$. By Corollary \ref{CPCor1}, statement (ii) holds. \\
 (ii) implies (iii):  It is immediate from Corollary \ref{CPCor1} that $\vp$ is completely positive.  For $u \in A_{\F^\perp} = (\pi(e_G) - z_F) \cdot A$ and $h \in Y$, 
 $$\vp(u)(h) = \l \alpha(h) , u \r = \l \alpha(h)  z_F, (\pi(e_G)-z_F) \c u \r = \l \alpha(h) z_F (\pi(e_G)-z_F), u \r = 0.$$ Hence, $\vp(u) = 0$.  \\
 (iii) implies (i): By Lemma   \ref{PositiveHomLem}, $e_H \in Y$ and, since $\vp$ is positive,  $\alpha(e_H) = \vp^*(e_H) \geq 0$; hence, $\alpha(e_H) z_F = z_F$ by Corollary \ref{Walter min idempotent etc Corollary}. Taking any $u \in A$, 
 \beqs \l \alpha(e_H) , u \r & = & \l \alpha(e_H),  (\pi(e_G)- z_F) \c u + z_F \c u \r \\
 & = & \vp( (\pi(e_G) - z_F) \c u)(e_H) +   \l \alpha(e_H) z_F,  u \r   =  \l z_F, u\r. \eeqs  Thus, $\alpha(e_H) = z_F$.   \end{proof} 
 
In \cite{Wal2}, critical idempotents $p$ in $\DBG$ are introduced and characterized. For instance, $p$ is critical precisely when $\bD_p$ is open in $p \DBG p$.   Thus, if $\vp:B(G) \ra B(H)$ is a completely positive  homomorphism and $\alpha(e_H) \in \bD_p$ for some critical idempotent $p$,  then $H_p = \alpha^{-1}(\bD_p)$ is an open (and closed) subgroup of $H_\alpha$---and therefore of $H$---by Theorem \ref{CPMainThm}.  We thus obtain the following immediate corollary to Theorem \ref{CPMainThm}; it, and several of the other results that follow, is also valid for 4-positive homomorphisms $\vp=\ja$ for which $\alpha(H_\alpha)$ is contained in $p\DBG p$.   Examples of critical idempotents in $\DBG$ are $e_G$ and $z_F$. 
 
\bc \label{CPConnectedCor1} Suppose that $\vp= \ja: B(G) \ra B(H)$ is a completely positive homomorphism and $\alpha(e_H) = \vp^*(e_H) \in \bD_p$ for some critical idempotent $p$. If $H_\alpha$ is connected---e.g., $H_\alpha = H$ is connected when $H$ is connected---then $\alpha:H_\alpha \ra \bD_p$ is a continuous homomorphism.  
\ec

 It might be tempting to conjecture that  $\alpha(Y)$ must always be contained in a subgroup of $\Delta(A)$ (and $\alpha$ is a homomorphism) when $\vp=\ja$ is completely positive.  (Of course, this condition is automatically satisfied when dealing with homomorphisms $\ja:A(G) \ra B(H)$.)  Perhaps surprisingly, as a corollary to one of our main theorems, we shall see that this is only true  when $G$ is compact (Corollary \ref{Compactness description corollary}).   We proceed in this direction  with some observations regarding the interplay between the different subgroups $\bD_p$ of $\Delta(A)$. 
 
 Let ${\cal P}_{\Delta(A)}$ denote the projections, (equivalently idempotents), in $\Delta(A)$.  Note that for $p, q \in \PDA$, $pq\in \PDA$ exactly when $pq=qp$ and $q\leq p$ if and only if $pq=qp = q$.  For a subset $D$  of $\Delta(A)$,  $D'=\{x \in \Delta(A) : xd = dx \ {\rm for \ all \ } d \in D\}$ is the commutant of $D$ in $\Delta(A)$.  
 
 \blem \label{Subgroup Interaction Prop} Suppose that $p, q \in \PDA$ and $pq = qp$. The following statements hold:  
 \bi \item[(i)] $pq \in \bD_p'\cap\bD_q'$ if and only if $p \in \bD_q'$ and $q \in \bD_p'$. 
 \item[(ii)] When (i) holds,  $\bD_p \bD_q, \bD_q \bD_p \subseteq \bD_{pq}$. 
 \item[(iii)] If $q\leq p$ and $q \in \bD_p'$, then $\bD_p \bD_q = \bD_q \bD_p = \bD_q$. 
 \ei 
 \elem 
 \begin{proof}  This is readily established. 
 \end{proof} 

Recall that $S$ is called a Clifford semigroup if each element in the semigroup $S$ belongs to a subgroup of $S$ and that a commutative idempotent semigroup is called a  semilattice. 
  
  \bp \label{Delta(A) as graded Clifford semigroup} If $\Delta(A)$ is a Clifford semigroup, then $$\Delta(A) = \dot{\bigcup}_{p \in \PDA} \bD_p.$$ Moreover, if  $\PDA$ is contained in $Z(\Delta(A))$, the algebraic centre of $\Delta(A)$, then  $\Delta(A)$ ``graded'' over the semilattice $\PDA$ in the sense that for each $p,q \in \PDA$,   $\bD_p \bD_q, \bD_q \bD_p \subseteq \bD_{pq}$.
  \ep 
  
  \begin{proof}  This is an immediate consequence of Proposition \ref{SpectralSubgroupsProp} and Lemma \ref{Subgroup Interaction Prop}.   \end{proof} 
  
  \br  \label{Delta(A) as graded Clifford semigroup Remark}   \rm 1.  In \cite{Il-Sp2}, the authors  showed that the spine compactification of $G$,    $G^*= \Delta(A)$ where  $A= A^*(G)$ is the spine of $B(G)$, is a Clifford semigroup with central idempotents; moreover, they provide a decomposition of $G^*$ as a union of groups $G_\tau$  graded over a semilattice  ${\cal T}_{nq}(G)$. Proposition \ref{Delta(A) as graded Clifford semigroup} thus gives an alternative description of the decomposition of $G^*$ into the subgroups $\bD_p$ of $ \Delta(A) \subseteq {\cal L}(\H_{\lambda^*})$  over the semilattice of projections $\PDA$ with its usual ordering.  Moreover, letting $\chi: G^*\ra \Delta(A)$ be the topological isomorphism described on page 286 of \cite{Il-Sp2}, for any $\tau \in {\cal T}_{nq}(G)$, $\chi(e_\tau)$ is an idempotent in $\Delta(A)$, so $\chi(e_\tau) = p_\tau$ for some unique $p_\tau$ in $\PDA$; since $\chi(G_\tau)$ is a subgroup of $\Delta(A)$ with identity $p_\tau$, $\chi(G_\tau) \subseteq \bD_{p_\tau}$ by Proposition \ref{SpectralSubgroupsProp}. As $\chi$ is one-to-one, $\chi(G_\tau) \cap \chi(G_{\tau'}) \subseteq \bD_{p_\tau} \cap \bD_{p_{\tau'}} = \emptyset$ for $\tau \neq \tau'$ and thus, since $\chi$ is also surjective, $\tau \mapsto p_\tau$ maps $ \tau \in {\cal T}_{nq}(G)$ onto $\PDA$ and $\chi$ maps each $G_\tau$ onto $\bD_{p_\tau}$ as a topological isomorphism. For $\tau, \tau' \in {\cal T}_{nq}(G)$, $p_\tau p_{\tau'} = \chi(e_\tau e_{\tau'}) \in \chi(G_\tau G_{\tau'}) \subseteq \chi(G_{\tau \wedge \tau'}) = \bD_{p_{\tau \wedge \tau'}}$, so $p_\tau p_{\tau'} = p_{\tau \wedge \tau'}$; thus ${\cal T}_{nq}(G)$ and $\PDA$ are isomorphic as semilattices.\\
  2.  Corollary  \ref{Delta(A) = G cupdot Gap Cor} and Example \ref{AFG Type Examples} provide examples of this decomposition of $\Delta(A)$. \\
  3.   Theorem 3.1 of  \cite{Liu-Mis} describes  a class of analytic groups $G$ for which the above decomposition of $\Delta(B(G))$ holds.\\
  4. Of course, $\PDA$ is contained in $Z(\Delta(A))$  whenever $G$ is abelian. 
  \er
  
  In the next lemma and the proposition that follows, we assume that $\vp=j_\alpha: A \ra B(H)$ is a completely positive homomorphism; thus, $\alpha:H_\alpha \ra \Delta(A)$ is a continuous map defined on an open subgroup $H_\alpha$ of $H$. We  further suppose that $\alpha(e_H) \in\bD_p$, so that $H_p = \{ k \in H_\alpha : \alpha(k) \in \bD_p\}$ is a subgroup of $H_\alpha$ and the conclusions of Theorem  \ref{CPMainThm}, which shall be used without further explanation, hold. (One can instead assume  that $\vp$ is 4-positive and $\alpha(H_\alpha) \subseteq p \Delta(A)p$, and replace conclusions of complete positivity with 4-positivity.)

  \blem \label{CP Complementary Thm Lemma} Suppose that  $y \in H_\alpha$ and $\alpha(y) \in \bD_q$ for some $q \in \PDA$. Then: 
  \bi \item[(i)] $q \leq p$,   $\alpha(yH_p) \subseteq \bD_q^{\ell} = \{ x \in \Delta(A): x x^* =q\}$  and \\   $\alpha(H_p y) \subseteq \bD_q^r = \{ x \in \Delta(A):  x^*x =q\}$;
  \item[(ii)] $\alpha(yH_p  \cup H_p y) \subseteq \bD_q$ if either (a) $yH_p = H_py$ or (b) $q \in \alpha(H_p)'$.     
  \ei  
  \elem
  
  \br   \label{CP Complementary Thm Lemma Remarks} \rm  Observe that (ii)(a) holds if $y$ is central in  $H_\alpha$ or $H_p$ is normal in $H_\alpha$---e.g., if $H_\alpha$ is abelian---and (ii)(b) holds if $q$ is central in $\Delta(A)$---e.g, if $A=A^*(G)$ or for any $A$ when $G$ is abelian. 
  \er
  
  \begin{proof} We have $pq= \alpha(e_H) \alpha(y)\alpha(y)^* = \alpha(e_Hy)\alpha(y)^* =q$ and for any $k \in H_p$, 
  $$\alpha(yk)\alpha(yk)^* = \alpha(y) \alpha(k)\alpha(k)^* \alpha(y)^* = \alpha(y) \alpha(kk^{-1})\alpha(y)^* = \alpha(ye_H)\alpha(y)^* = q;$$
similarly, $\alpha(H_py) \subseteq \bD_q^r$. We have established (i).  If (ii)(a) holds, $\alpha(yH_p) \subseteq \bD_q$ by (i) and if (ii)(b) holds, then for $k \in H_p$, 
$$\alpha(yk)^* \alpha(yk) = \alpha(k)^* \alpha(y)^* \alpha(y) \alpha(k) =  \alpha(k)^* q \alpha(k) =  \alpha(k)^*  \alpha(k) q = pq = q. \qedhere$$ 
  \end{proof} 
   
   For each $q \in \PDA$ for which $\alpha^{-1}(\bD_q)$ is nonempty, let $$\alpha_q : H_\alpha \ra \Delta(A): h \mapsto q \alpha(h) q  \quad {\rm and} \quad H_q = \alpha_q^{-1}(\bD_q).$$
   Observe that since $\alpha(H_\alpha)$ is contained in $p \Delta(A)p$, $\alpha_p = \alpha$ and our definition of $H_q$ is consistent with the definition of $H_p$ introduced in the statement of Theorem \ref{CPMainThm}.  
   
   The following proposition provides more general information about the behaviour of completely positive homomorphisms and will  be used in the proof of Theorem \ref{Second Main Thm} below.  Lemma \ref{CP Complementary Thm Lemma} and Remarks \ref{CP Complementary Thm Lemma Remarks} provide sufficient conditions for when the assumption in statement (iii)  holds. 
      
 \bp \label{CP Complementary Prop} Suppose that $\alpha^{-1}(\bD_q)$ is nonempty. Then: 
 \bi \item[(i)] $q \leq p$, $j_{\alpha_q} : A \ra B(H)$ is a completely positive homomorphism, $H_q$ is a subgroup of $H_\alpha$ such that $\alpha_q : H_q \ra \bD_q$ is a homomorphism and for each $h \in H_\alpha, \, k \in H_q$, 
 $$\alpha_q (hk) = \alpha_q(h) \alpha_q(k), \quad \alpha_q(kh) = \alpha_q(k) \alpha_q(h) \ \ and  \ \ \alpha_q(h^{-1}) = \alpha_q(h)^*  . $$ 
     \item[(ii)]  $\alpha^{-1}( \bD_q) \subseteq H_q$  and $\alpha$ agrees with $\alpha_q$ on $\alpha^{-1}(\bD_q)$; thus $\alpha_q{\large |}_{H_q}$ is a homomorphic extension of $\alpha{\large |}_{\alpha^{-1}( \bD_q)}$ to the subgroup $H_q$. 
 \item[(iii)] If $\alpha(y_q H_p) \subseteq \bD_q$ for some $y_q \in H_\alpha$, then  $H_p  \subseteq H_q$ and for every  $y$ such that $\alpha(y)\in \bD_q$, $\alpha(yH_p \cup H_p y ) \subseteq \bD_q$. 
 \ei 
 \ep   
 
 \begin{proof}   (i) By Lemma \ref{CCHomLem}, $j_{\alpha_q}$ is a completely contractive homomorphism and by Lemma \ref{CP Complementary Thm Lemma}, $q \leq p$, so $\alpha_q(H_\alpha) = q \alpha(H_\alpha)q  \subseteq qp \Delta(A)pq \subseteq q VN_\pi q$, where $A = A_\pi$.  As with the proof of Theorem \ref{CPMainThm} under hypothesis (ii), it follows that $(j_{\alpha_q})^*$ maps $W^*(H)$ into the von Neumann subalgebra $q VN_\pi q$ of $VN_\pi$ and $(j_{\alpha_q})^*(e_H) = \alpha_q(e_H) = q$, the identity element of $q VN_\pi q$.  By Theorem \ref{CPBasicThm}(b),  ($j_{\alpha_q})^*$, and therefore $j_{\alpha_q}$, is a completely positive homomorphism.  
 The rest of (i) follows immediately from Theorem \ref{CPMainThm}. 
 
 (ii) If $\alpha(h) \in \bD_q$, then $\alpha_q(h) = q \alpha(h)q = \alpha(h)$. 
 
 (iii) Suppose that $\alpha(y_qH_p) \subseteq \bD_q$.   Then for $k \in H_p$, $$\alpha_q(k) = q \alpha(k) q = \alpha(y_q)^* \alpha(y_q) \alpha(k)q = \alpha (y_q)^* \alpha(y_qk) q \in \bD_q,$$ so $k \in H_q$.  Also, $$\alpha(ky_q) = \alpha(k) \alpha(y_q) = \alpha(k)q\alpha(y_q) =  \alpha(k)\alpha(y_q)^*\alpha(y_q)^2= \alpha(y_q k^{-1})^* \alpha(y_q)^2 \in \bD_q,$$ so $\alpha(H_p y_q)\subseteq \bD_q$.  Now for any $y\in H_\alpha$ such that $\alpha(y) \in \bD_q$, 
 $$\alpha(yk) = \alpha(y)\alpha(k) = \alpha(y)q \alpha(k) = \alpha(y) \alpha(y_q)^*\alpha(y_q)\alpha(k) = \alpha(y) \alpha(y_q)^*\alpha(y_q k),$$  which belongs to  $\bD_q$. Since we have shown that $\alpha(H_p y_q) \subseteq  \bD_q$, we similarly obtain $\alpha(H_p y ) \subseteq \bD_q$. 
 \end{proof}

 The next two lemmas will be used in the proof of Theorem \ref{Second Main Thm}.
 
 \blem  \label{Lemma for Second Main Thm} Let $p_1, \ldots, p_n \in \PDA$, $$  L = \bigcupdot_{i=1}^n \bD_{p_i} \quad {\rm and} \quad   L_k = \bigcupdot_{i=k}^n \bD_{p_i}.$$   Suppose that $p_1 \geq p_2 \geq \ldots \geq p_n$ and $p_{i+1} \in \bD_{p_i}'$ for $i=1, \ldots, n-1$.  Then  each $p_i $ is central in $L$,   $L$ is a semigroup with identity $p_1$,  and  $ L_k $ is an ideal in $L$.  
 \elem
 
 \begin{proof}  Let $x \in L$, say $x \in \bD_{p_i}$.  Since $p_1\geq p_i$, $p_1x = p_1(p_ix) = (p_1p_i)x = p_ix = x$ and similarly, $x p_1 = x$. In particular, $p_1 \in L'$. Suppose that $p_1, \ldots, p_k \in L'$. If $i \geq k+1$, then---as above---$p_{k+1} x = x = xp_{k+1}$ since $p_{k+1} \geq p_i$, so suppose $i \leq k$. By hypothesis, $p_k \in \bD_{p_i}'$ and $p_i \geq p_k$, so $xp_k = p_kx \in \bD_{p_k}$ by Lemma \ref{Subgroup Interaction Prop}(iii). Since $p_k \geq p_{k+1}$ and $p_{k+1} \in \bD_{p_k}'$, we therefore obtain
 $$p_{k+1}x = (p_{k+1}p_k)x = p_{k+1}(p_k x)  =  (p_k x) p_{k+1} =  (x p_k ) p_{k+1} =  x (p_k  p_{k+1}) = x p_{k+1}. $$ Thus, each $p_i$ is central in $L$.  Hence, $p_j \in \bD_{p_i}'$ for every $i$ and $j$, so $ \bD_{p_i} \bD_{p_j} =  \bD_{p_{i \vee j}}$ by Lemma   \ref{Subgroup Interaction Prop}; the remaining statements quickly follow. 
 \end{proof} 
  
  \blem  \label{Discrete CP implies CP Lemma} Let $H_d$ denote the group $H$ with the discrete topology, $\iota_d: B(H) \hookrightarrow B(H_d)$,   $\vp: A\ra B(H)$ a linear map. If  $\vp_d = \iota_d \circ \vp: A \ra B(H_d)$ is  completely positive, then  $\vp$ is completely positive. 
  \elem
  
  \begin{proof}  Since $B(H)$ is closed in $B(H_d)$ \cite[(2.24)]{Eym}, there is a representation  $(\omega_H)^\perp$ of $H_d$, disjoint from $\omega_H$ (viewed here as a representation of $H_d$), such that $B(H_d) = A_{\omega_H} \oplus_1 A_{(\omega_H)^\perp}$ and $W^*(H_d) = VN_{\omega_H} \oplus_{\infty} VN_{(\omega_H)^\perp}$ \cite[3.12, 3.13, 3.18]{Ars}. The dual map of $\iota_d$ is the $*$-homomorphism $\iota_d^*:W^*(H_d) \ra W^*(H): x\oplus z \mapsto x$, so for $x_1, \ldots, x_n \in W^*(H)$ and $y_1, \ldots, y_n \in A^* = VN_\pi$, 
  $$\sum_{i,j=1}^n y_i^*\vp^*(x_i^* x_j)y_j = \sum_{i,j=1}^n y_i^*(\vp_d)^*((x_i\oplus 0)^* (x_j\oplus0))y_j.$$
 By \cite[Corollary IV.3.4]{Tak}, $\vp$ is completely positive when $\vp_d$ is completely positive. 
  \end{proof}  
  
  We introduce some new terminology to help with the phrasing of the next several results:
  
  \bd \rm Let $p_1, \ldots, p_n \in \PDA$ be such that  $p_1 \geq p_2 \geq \ldots \geq p_n$ and $p_{i+1} \in \bD_{p_i}'$ for $i=1, \ldots, n-1$.  By a \it compatible sytem  of homomorphisms (affine maps) for $p_1, \ldots, p_n$, \rm we shall mean a sequence 
  $$ \beta_i: H_i \ra \bD_{p_i} \qquad ( \beta_i: E_i \ra \bD_{p_i})$$
  of homomorphisms (affine maps) defined on subgroups $H_i$ (cosets $E_i$) of $H$ such that 
   $$H_1 <  H_2 <  \ldots <  H_n   \qquad ( E_1 \subset E_2 \subset  \ldots \subset E_n)$$ and for each $h \in H_i$ ($h \in E_i$), 
  $$\beta_{i+1}(h) = p_{i+1} \beta_i(h) =  p_{i+1} \beta_i(h) p_{i+1}.$$  The map $$\alpha = \gamma_1 \cupdot \gamma_2 \cupdot \cdots \cupdot \gamma_n,$$  where  $\gamma_1 = \beta_1$   and   $\gamma_i = \beta_i {\large |}_{H_i\bs H_{i-1}}$ for $i = 2, \ldots, n$, will be called the \it fusion map \rm of $\beta_1, \ldots, \beta_n$. 
  \ed

  \bt  \label{Second Main Thm}  Let $p_1, \ldots, p_n \in \PDA$ be such that  $p_1 \geq p_2 \geq \ldots \geq p_n$ and $p_{i+1} \in \bD_{p_i}'$ for $i=1, \ldots, n-1$. If $\ds \alpha: Y \ra \bigcupdot_{i=1}^n \bD_{p_i}$ is such that each $\alpha^{-1} (\bD_{p_i})$ is nonempty, then the following statements are equivalent: 
  \bi  \item[(i)] $\vp = \ja$ is a completely positive homomorphism of $A$ into $B(H)$; 
  \item[(ii)] $\vp = \ja$ is a 4-positive homomorphism of $A$ into $B(H)$ and $\alpha(e_H) \in \bD_{p_1}$; 
  \item[(iii)]  $\alpha$ is continuous, $Y$ is an open subgroup of $H$, and $\alpha$ is the fusion map of a compatible system of homomorphisms for $p_1, \ldots, p_n$.  
  \ei 
  \et

  \begin{proof} (i) implies (ii): By Lemma \ref{PositiveHomLem}, $\alpha(e_H) = p_k$ for some $k$, and by Proposition \ref{CP Complementary Prop}, $p_i \leq p_k$ for each $i$. Hence, $p_k = p_1$.  
 
 \smallskip

 \noindent  (ii) implies (iii): By 
  Lemma \ref{PositiveHomLem}, $Y=H_\alpha$ is an open subgroup of $H$,  $\alpha$ is continuous on $H_\alpha$ and $\alpha(e_H) = p_1$. As observed in Lemma \ref{Lemma for Second Main Thm}, $p_1$ is the identity for the semigroup $L= \bigcupdot_1^n \bD_{p_i}$, so $\alpha(H_\alpha)$ is contained in $p_1\Delta(A)p_1$; hence all of the conclusions of Theorem \ref{CPMainThm} and Proposition \ref{CP Complementary Prop} hold.   For each $i$, letting $$\alpha_i: H_\alpha \ra L: h \mapsto p_i\alpha(h)$$
---equivalently, $\alpha_i(h) = p_i \alpha(h) p_i$, since each $p_i \in L'$---$H_i = \alpha_i^{-1} (\bD_{p_i})\leq H_\alpha$  and $\beta_i:= \alpha_i{\large |}_{H_i} :H_i \ra \bD_{p_i}$ is a homomorphism by Proposition \ref{CP Complementary Prop}. Observe that $H_n = H_\alpha$, since $L_n = \bD_{p_n}$ is an ideal in $L$, and $$\alpha_{i+1}(h) = p_{i+1}\alpha(h) = p_{i+1}p_i \alpha(h) = p_{i+1} \alpha_i(h) \qquad (h \in H_\alpha).$$ Hence, for $h \in H_i$, Lemma \ref{Subgroup Interaction Prop}(iii) gives $\alpha_{i+1}(h) \in \bD_{p_{i+1}} \bD_{p_i} = \bD_{p_{i+1}}$, and therefore $h \in H_{i+1}$; thus, $H_i \leq H_{i+1}$ and $\beta_{i+1} (h)  = p_{i+1} \beta_i(h)$ for $h \in H_i$. Let $h \in H_\alpha = H_n$. If $h \in H_1$, then $\alpha(h) = \alpha_1(h) = \gamma_1(h)$, so suppose that $h \in H_k \bs H_{k-1}$ for some $k \geq 2$. If $\alpha(h) \in \bD_{p_j}$ for  $j < k$, then $\alpha_j(h) = p_j \alpha(h) \in \bD_{p_j}$, giving $h \in H_j \leq H_{k-1} $, a contradiction; if $\alpha(h)$ belongs to the ideal $L_{k+1}$ of $L$, $\alpha_k(h) = p_k \alpha(h) \in L_{k+1}$, also a contradiction. It follows that $\alpha(h) \in \bD_{p_k}$ and we obtain $\alpha(h) = p_k\alpha(h) = \alpha_k(h) = \gamma_k(h)$, as needed.  

 \smallskip

 \noindent  (iii) implies (i): For notational purposes, denote the empty set by $H_0$.  For each $i$, we know that $j_{\beta_i}: A \ra B(H_d)$ is a (completely positive) homomorphism by Corollary \ref{CPHomCor}. Also, $1_{H_i \bs H_{i-1}}$ is an idempotent in $B(H_d)$, so ${\mathfrak 1}_i: B(H_d) \ra B(H_d): u \mapsto u 1_{H_i\bs H_{i-1}}$ is an automorphism of $B(H_d)$, and therefore $$j_{\gamma_i} = \mathfrak{1}_i \circ j_{\beta_i}: A \ra B(H_d) \qquad (i = 1, \ldots, n) $$  is a disjoint sequence of homomorphisms. By Lemma \ref{Disjoint Homom Lemma}, $\sum j_{\gamma_i} = j_{\cupdot \gamma_i} = j_\alpha$ is a homomorphism of $A$ into $B(H_d)$. However, $H_n=Y$ is an open (and therefore closed) subgroup of $H$, and $\alpha: H_n \ra \Delta(A)$ is continuous, so it readily follows that $j_\alpha(u)$ is continuous on $H$ for any $u \in A$; therefore $j_\alpha(u)$ belongs to $B(H)$ by \cite[(2.24)]{Eym}. Thus, $\vp =j_\alpha: A \ra B(H)$ is a homomorphism. 
 
 By Lemma \ref{Discrete CP implies CP Lemma}, it now suffices for us to show that $\vp_d = j_\alpha: A \ra B(H_d)$ is completely positive. To this end, we \it claim \rm that $\vp_d^*: W^*(H_d) \ra A^* =VN_\pi$ satisfies 
 \beq \label{Second Main Thm Eqn}  \vp_d^*(x) = \sum_{i=1}^n (p_i - p_{i+1}) j_{\beta_i}^*(x)(p_i-p_{i+1}) \qquad (x \in W^*(H_d)),
 \eeq
 where $p_{n+1} = 0$.  By linearity,  wk$^*$-continuity and wk$^*$-density of $\l H_d\r$ in $W^*(H_d)$, it suffices to establish (\ref{Second Main Thm Eqn}) for $h \in H_d$. Obviously both sides of (\ref{Second Main Thm Eqn}) are zero when $h \in H\bs H_n$, so suppose that $h \in H_k\bs H_{k-1}$ for some $k$. Then 
 \beqs \sum_{i=1}^n (p_i - p_{i+1}) j_{\beta_i}^*(h)(p_i-p_{i+1}) & = & \sum_{i=k}^n (p_i - p_{i+1}) \beta_i(h)(p_i-p_{i+1}) \\
 & = &   \sum_{i=k}^{n-1} \left( p_i\beta_i(h)  - p_{i+1} \beta_i(h)\right) + p_n \beta_n(h)\\
 & = & \sum_{i=k}^{n-1} \left( \beta_i(h)  -  \beta_{i+1}(h)\right) +  \beta_n(h)  \\
 & = & \beta_k(h) = \alpha(h) = \vp_d^*(h),
 \eeqs 
 as required. Since each $j_{\beta_i}^*$ is completely positive, conjugation by an operator is completely positive, and each $p_i - p_{i+1}$ is a projection, it follows from (\ref{Second Main Thm Eqn}) that $\vp_d^*$ is completely positive. 
  \end{proof}

 
 \bex    \label{Construction of compatible homoms Ex}  \rm    Let $p_1, \ldots, p_n \in \PDA$ be such that  $p_1 \geq p_2 \geq \ldots \geq p_n$ and $p_{i+1} \in \bD_{p_i}'$ for $i=1, \ldots, n-1$.  Let $H$ be any group containing a sequence of  proper subgroups 
 $H_1 < H_2 < \cdots < H_n$, and let $\beta: H_n \ra \bD_{p_1}$ be any homomorphism. For each $i = 1, \ldots, n$,  define $\beta_i:H_i \ra \bD_{p_i}$ by putting $\beta_i(h) = p_i \beta(h)$. Since $\bD_{p_i}$ is an ideal in $\bD_{p_1} \cupdot \bD_{p_i}$ and $p_i \in \bD_{p_1}'$ (Lemma \ref{Lemma for Second Main Thm}),  $\beta_i$ is a homomorphism of $H_i$  into  $\bD_{p_i}$ and,  since $p_{i+1} \leq p_i$, $\beta_{i+1}(h) = p_{i+1} \beta_i(h)$ for $h \in H_i$.   Hence, $\beta_1, \ldots, \beta_n$  is a compatible system of homomorphisms for $p_1, \ldots, p_n$. If, further, $H_1, \ldots, H_n$ are open in $H$ (they are trivially so if $H$ is discrete) and we choose $\beta$ to be continuous (e.g., the trivial homomorphism), then the fusion map $\gamma$ is clearly continuous. Thus, by (iii) implies (i) of Theorem \ref{Second Main Thm}, it is easy to construct completely positive homomorphisms $\ja : A \ra B(H)$ where  $\alpha$ does   not  map into a subgroup of $\Delta(A)$. 
  \eex
 
 The following corollary improves on \cite[Corollary 5.2]{Il-St}. 
 
 \bc \label{Compactness description corollary}  Let $G$ be a locally compact group. For any disconnected locally compact group $H$, the following statements are equivalent:  \bi \item[(i)] $G$ is compact;
 \item[(ii)] for any completely bounded homomorphism $\ja:B(G) \ra B(H)$, $\alpha$ maps into a subgroup of $\DBG$; 
 \item[(iii)] for any completely positive homomorphism $\ja:B(G) \ra B(H)$, $\alpha$ maps into a subgroup of $\DBG$.
 \ei
 \ec 
 
 \begin{proof}   If $G$ is compact, $\DBG = G$, so (i) implies (ii) is trivial, as is (ii) implies (iii).  If $G$ is non-compact, then, by Corollary \ref{Walter min idempotent etc Corollary}, $e_G$ and $z_F$ are distinct central idempotents in $\DBG$ with $e_G \geq z_F$. Let $H_1$ be a proper open subgroup of $H$, $\beta: H \ra G$ the trivial homomorphism, $\beta_1 = \beta \large{|}_{H_1} :H_1 \ra G$, $\beta_2: H \ra \bD_{z_F}: h \mapsto z_F \beta(h) = z_F$. By Example   \ref{Construction of compatible homoms Ex}, if $\alpha$ is the fusion map of $\beta_1$ and  $ \beta_2$, then $\ja:B(G) \ra B(H)$ is a completely positive homomorphism such that $\alpha$ maps into both $G$ and $\bD_{z_F}$. 
 \end{proof}

 Under the hypothesis of Theorem \ref{Second Main Thm}, we can  also characterize the completely contractive homomorphisms $\vp: A \ra B(H)$:
 
 \bt  \label{CCThm} Let $p_1, \ldots, p_n \in \PDA$ be such that  $p_1 \geq p_2 \geq \ldots \geq p_n$ and $p_{i+1} \in \bD_{p_i}'$ for $i=1, \ldots, n-1$. If $\ds \alpha: Y \ra \bigcupdot_{i=1}^n \bD_{p_i}$ is such that each $\alpha^{-1} (\bD_{p_i})$ is nonempty, then the following statements are equivalent: 
  \bi  \item[(i)] $\vp = \ja$ is a completely contractive homomorphism of $A$ into $B(H)$; 
  \item[(ii)]  $\alpha$ is continuous, $Y$ is an open coset in $H$, and  there exists an element $y_0 \in \alpha^{-1}(\bD_{p_1})$ and $x_0 \in \bD_{p_1}$ such that 
  $$\alpha(y) = x_0 \gamma(y_0^{-1}y) \qquad (y \in Y)$$ 
  where  $\gamma$ is the fusion map of a compatible system of homomorphisms for $p_1, \ldots, p_n$.  
  \item[(iii)]  $\alpha$ is continuous, $Y$ is an open coset in $H$, and $\alpha$ is the fusion map of a compatible system of affine maps  for $p_1, \ldots, p_n$.  
  \ei 
  \et 

\begin{proof} (i) implies (ii):   Assuming (i), $1_Y = \vp(1_G)$  is a norm-one idempotent, so $Y$ is an open coset in $H$ by Corollary \ref{Il-SpThm2.1Cor}. Take $y_0 \in \alpha^{-1}(\bD_{p_1})$ and define $\gamma$ on the open  subgroup $H_\alpha:= y_0^{-1}Y$ of $H$ by
$$\gamma: H_\alpha \ra  L= \bigcupdot_{i=1}^n \bD_{p_i}: h \mapsto \alpha(y_0)^* \alpha(y_0h). $$
By Lemma \ref{CCHomLem},  $j_\gamma: A \ra B(H)$ is a completely contractive homomorphism, and $j_\gamma^*(e_H) = \gamma(e_H) = \alpha(y_0)^* \alpha(y_0) = p_1$. Moreover, $j_\gamma^*(H)$ is contained in $L \cup \{0\}= p_1(L \cup \{0\})p_1$, so $j_\gamma^*$ is a completely contractive identity-preserving linear mapping of $W^*(H)$ into the von Neumann algebra $p_1 VN_\pi p_1$. By Theorem \ref{CPBasicThm}, $j_\gamma^*$ is completely positive, so by Theorem \ref{Second Main Thm}, $\gamma$ is the fusion map of a compatible system of homomorphisms for $p_1, \ldots, p_n$. Finally, for any $y = y_0 h \in y_0  H_\alpha= Y$, 
$$\alpha(y_0) \gamma(y_0^{-1}y) = \alpha(y_0) \alpha(y_0)^* \alpha(y_0h) = p_1 \alpha(y) = \alpha(y).  $$  
\smallskip 

\noindent (ii) implies (iii): Suppose that $\gamma$ is the fusion map of homomorphisms $\beta_i: H_i \ra \bD_{p_i}$ ($i=1, \ldots, n$). For each $i$, let $E_i = y_0H_i$, 
$$\alpha_i: E_i \ra \bD_{p_i}: y \mapsto x_0 \beta_i(y_0^{-1}y).$$
Equivalently, $\alpha_i(y) = x_0 p_i \beta_i (y_0^{-1}y)$ and, since $x_0 \in \bD_{p_1}$,  $x_0 p_i \in \bD_{p_i}$; we  conclude that  each $\alpha_i$ is affine.  The $p_i$ are central in $L$, so $\alpha_1, \ldots, \alpha_n$ is a compatible system of affine maps for $p_1, \ldots, p_n$, and it is easy to see that $\alpha$ is the associated fusion map. 

\smallskip 

\noindent (iii) implies (ii): Suppose that $\alpha$ is the fusion map of  $\alpha_i : E_i \ra \bD_{p_i}$ ($i=1, \ldots, n$), choose $y_0 \in E_1$ and let $x_0 = \alpha(y_0)$.  For each $i$, let $H_i = y_0^{-1} E_i$, $$\beta_i: H_i \ra \bD_{p_i}: h \mapsto x_0^* \alpha_i(y_0h).$$ Then it is easy to see that $\beta_1, \ldots, \beta_n$ is a compatible system of homomorphisms for $p_1, \ldots, p_n$ and,  letting $\gamma$ denote the associated fusion map, one can also readily check that $\alpha(y) = x_0\gamma(y_0^{-1} y)$ for each $y \in E_n=Y$. 

\smallskip 

\noindent (ii) implies (i):  It follows from (ii) that $\gamma$ is defined on the open subgroup $H_\alpha = y_0^{-1} Y$ of $H$ and $\gamma(h) = x_0^* \alpha(y_0h)$, so $\gamma$ is continuous. By Theorem \ref{Second Main Thm}, $j_\gamma: A \ra B(H)$ is a completely positive (and therefore completely contractive) homomorphism.  By Lemma \ref{CCHomLem}, we conclude that $\vp = \ja = \ell_{y_0^{-1}} \circ j_\gamma \circ \ell_{x_0}$ is a completely contractive homomorphism of $A$ into $B(H)$.    \end{proof}

 \br \rm Under the hypothesis of Theorem \ref{CCThm}, we showed that $\vp:A \ra B(H)$ is a completely contractive homomorphism exactly when it has a Cohen-type factorization $\vp = \ell_{y_0} \circ j_\gamma \circ \ell_{x_0}$ where $\gamma$ is the fusion map of a compatible system of homomorphisms, $y_0 \in H$ and $x_0 \in \Delta(A)$; cf. the main results concerning (completely) contractive homomorphisms in  \cite{Gli}, \cite{Il-Sp} and \cite{Pham}.  
 As well, we have the following general statement, which is a  consequence  of Proposition \ref{SpectralSubgroupsProp} and Theorem \ref{CCThm}. 
 \er

\bc \label{CC Thm Cor}  Let   $\vp = \ja:A \ra B(H)$ be a homomorphism. The following statements are equivalent: \bi 
\item[(i)]   $\vp$ is completely contractive  and $\alpha(Y)$ is contained in a subgroup of $\Delta(A)$; 
\item[(ii)] there is an open subgroup $H_0$ of $H$, a continuous homomorphism $\gamma:H_0 \ra \Delta(A)$, and elements $x_0 \in \bD_{\gamma(e_H)}$ and $y_0 \in H$ such that $Y = y_0H_0$ and $\alpha(y) = x_0\gamma(y_0^{-1} y)$, $(y \in Y)$;  
\item[(iii)] $Y$ is an open coset in $H$ and   $\alpha: Y \ra \Delta(A)$ is a continuous affine mapping of $Y$  into some subgroup of  $\Delta(A)$.
\ei  
\ec

For several groups $G$ and subalgebras $A$ of $B(G)$, the hypotheses of Theorems \ref{Second Main Thm} and \ref{CCThm} are always satisfied and we thus obtain the full characterization of \it all  \rm completely positive and completely contractive homomorphisms in these situations. 

\bc  \label{Full characterization of cc and cp in special cases Cor}  Let $A$ be any closed, translation-invariant unital subalgebra of $B(G)$ for which the projections in $\Delta(A)$, $\PDA$, are  central in $\Delta(A)$ and linearly ordered, $\PDA$ is finite,  and  $$\Delta(A) = \bigcupdot_{p \in \PDA} \bD_{p}.$$  Then the following statements are equivalent:  
 \bi  \item[(i)] $\vp = \ja: A \ra B(H)$ is a  completely positive (completely contractive) homomorphism; 
  \item[(ii)]  $\alpha$ is continuous, $Y$ is an open subgroup (coset) in $H$, and $\alpha$ is the fusion map of a compatible system of homomorphisms (affine maps)   for some subset of  $\{p_1, \ldots, p_n\}$.  
  \ei 
Some cases in which the hypotheses  are satisfied are: 
\bi \item $A= B(G)$ when $G$ is either the Euclidean motion group $\mathbb{R}^n \rtimes SO(n) $ for some positive integer $n\geq 2$, or  $G$  is the $p$-adic motion group $\Q_p^n \rtimes SL(n, \mathbb{O}_p)$ for some  positive integer $n$ and some prime $p$;
\item $A= A^*(G)$, the spine of $B(G)$, when $G$ contains a compact normal subgroup $K$ for which $G/K$ is topologically isomorphic to  $\mathbb{R}$ or $\mathbb{Z}$;  the p-adic field $\mathbb{Q}_p$ for some prime number $p$;    any minimally weakly almost periodic group such as $G= SL_2(\mathbb{R})$ or $G=\mathbb{R}^n \rtimes SO(n)$; the $ax+b$-group. 
\item $A = \AFG$ for any locally compact group $G$. 
\ei 
\ec 
 
 \begin{proof} The equivalence is an immediate consequence of Theorems \ref{Second Main Thm} and \ref{CCThm}, so we only need to  observe that the hypotheses of the corollary are  satisfied in the specified cases:  As noted in Example \ref{AFG Type Examples}, except for when $A=A^*(G)$ and $G$ is the $ax+b$-group, in each case  $\Delta(A) = G \cupdot G^{ap} =   \bD_{e_G} \cupdot \bD_{z_F}$ and $e_G \geq z_F$. When $A = A^*(G)$, we observed in Remark \ref{Delta(A) as graded Clifford semigroup Remark}  that $A$ takes the form described in Proposition  \ref{Delta(A) as graded Clifford semigroup} with $\PDA$ central in $\Delta(A)$ and $\PDA \simeq {\cal T}_{nq}(G)$ as semilattices.  In the case that  $G$ is the $ax+b$-group, ${\cal P}_{\Delta(A^*(G))} = \{e_G, p, z_F\}$ with $e_G \geq p \geq z_F$ by \cite[Theorem 6.6]{Il-Sp2}.   
 \end{proof}

\br  \label{CP and CC Homom G cup Gap Remarks}   \rm  1.  Let $A$ be any of the examples from Corollary \ref{Full characterization of cc and cp in special cases Cor} for which $\Delta(A) = G \cupdot G^{ap}$, and let  $\vp=\ja: A \ra B(H)$ be a completely positive homomorphism. Then  $Y= H_\alpha$ is an open subgroup of $H$,  $\alpha:H_\alpha \ra G \cupdot G^{ap}$ and, by Corollary \ref{Full characterization of cc and cp in special cases Cor},   exactly one of the following three statements  holds: 
\bi  \item[(a)] $\alpha: H_\alpha \ra G^{ap}$ is a continuous homomorphism, (by Corollary \ref{CP Cor z_F}, this happens if and only if $\alpha(e_H) \in G^{ap}$); 
\item[(b)]  $\alpha: H_\alpha  \ra G$ is a continuous homomorphism; 
\item[(c)]  there is a proper subgroup $H_a$ of $H_\alpha$ and continuous homomorphisms
$$\alpha_a : H_a \ra G \quad {\rm and}  \quad \beta: H_\alpha \ra G^{ap}$$ 
such that for $h \in H_a$,  $\beta(h) = \alpha_a(h) e_{ap}$---equivalently, $\beta{\large |}_{H_a} = \delta^{ap} \circ \alpha_a$ by Corollary \ref{Walter min idempotent etc Corollary}(ii)---and $\alpha = \alpha_a \cupdot \alpha_s $ where $\alpha_s = \beta{\large |}_{H_\alpha\bs H_a}$.
\ei
Conversely, if $\alpha$ takes any of the forms described in statements (a), (b) and (c), then $\ja: A \ra B(H)$ is a completely positive homomorphism. 

One can similarly provide explicit descriptions of all possible formulations of $\alpha$ when $\vp$ is a completely contractive homomorphism.   When $G$ is the $ax+b$-group and $\ja:A^*(G) \ra B(H)$ is a completely positive/completely contractive homomorphism, explicit descriptions of all possible formulations of  $\alpha$ can also be readily formulated from Corollary \ref{Full characterization of cc and cp in special cases Cor}.    

\smallskip 

\noindent 2. It can  happen that  $\PDA$ is not linearly ordered, even when $\Delta(A)$ takes the form described in  Proposition \ref{Delta(A) as graded Clifford semigroup}. For example, this is the case when $A = A^*(\mathbb{R}^n)$ and $n\geq 2$ \cite[Theorem 6.3]{Il-Sp2}. 

\smallskip 

\noindent 3. Suppose that $\Delta(A)$ takes the form described in  Proposition \ref{Delta(A) as graded Clifford semigroup} and $\ja: A \ra B(H)$ is completely positive. If  $\alpha(e_H) \in \bD_{p}$ and we at least know that ${\cal P}_{\leq p} = \{ q \in \PDA: q \leq p\}$ is finite and linearly ordered, then by Proposition \ref{CP Complementary Prop} and Theorem \ref{Second Main Thm} $\alpha $ is a continuous fusion map of a  compatible system of homomorphisms   for some subset of   ${\cal P}_{\leq p}$.    
 \er 
 

 
 \section{Spatial homomorphisms revisited}

 Let $\vp = \ja:B(G) \ra B(H)$ be a homomorphism such that $\alpha(Y)$ is contained in $\bD_p$ where $p$ is a critical idempotent in $\DBG$. Then $\bD_p$ is a locally compact group \cite[Cor. to Proposition 7]{Wal2} and the embedding $\iota: \bD_p \hookrightarrow \DBG$ is a continuous homomorphism, so---by Proposition \ref{CPHomProp}---$$j_p \coloneqq j_\iota: B(G) \ra B(\bD_p), \qquad j_p u (x) = \l x, u \r_{W^*-B(G)}$$
 is a completely positive homomorphism and complete quotient mapping onto the closed, translation-invariant, point-separating, unital subalgebra
 $$A_{\omega_p} = \{ \xi *_{\omega_p} \eta: \xi, \eta \in \K = p \H_{\omega_G} \}$$ 
 of $B(\bD_p)$, where $\omega_p$ is the continuous unitary representation of $\bD_p$ defined by 
 $$\omega_p \coloneqq \omega_\iota : \bD_p \ra {\cal U}( \K): x \mapsto p_\K x E_\K (= p \, x\big{|}_\K).$$
 Observe that if we view $\K$ as a subspace of $\H_{\omega_G}$, each $u \in A_{\omega_p}$ can be written as 
 $$u(x) = \l \omega_p(x) \xi | \eta \r_\K = \l x \, \xi |\eta\r_{\H_{\omega_G}} \qquad (x \in \bD_p), $$
 and $\xi, \eta \in \K$ can be chosen so that $\| u\| = \| \xi \| \| \eta\|$. Letting 
 $$\alpha_p : Y \subseteq H \ra \bD_p : y \mapsto \alpha(y),$$
 $\alpha_p$ is a continuous mapping of the set $Y \in \Omega(H)$ into the locally compact group $\bD_p$.  
 
 \bp \label{SpatialHomFactorizationProp} Let $\vp = \ja:B(G) \ra B(H)$ be a homomorphism such that $\alpha(Y)$ is contained in $\bD_p$ where $p$ is a critical idempotent in $\DBG$. The map $j_{\alpha_p}: A_{\omega_p} \ra B(H)$ is a well-defined spatial homomorphism and we have the factorization $\ja = j_{\alpha_p} \circ j_p$; i.e., the  diagram 
  $$  \xymatrixrowsep{2pc} \xymatrixcolsep{3pc}
\xymatrix{ & & B(G)\ar@<.5ex>[d]_{j_p}  \ar@<.5ex>[rr]^{\vp = \ja}
 & & B(H)\\
& B(\bD_p) \supseteq  \mkern-72mu & A_{\omega_p}  \ar@<.5ex>[rru]_{j_{\alpha_p}} } $$
 commutes. Moreover, the spatial homomorphism $j_{\alpha_p}$ is positive/ contractive/ completely positive/ completely contractive/ completely bounded if and only if $\ja$ is so. 
 \ep 
 
 \br  \rm Thus, when $\alpha(Y)$ is contained in $\bD_p$ with $p$ critical, the problem of studying $\vp = \ja$ is reduced to the problem of studying the spatial homomorphism $j_{\alpha_p}: A_{\omega_p} \ra B(H)$. 
 \er
 
 \begin{proof}  Since $Y$ is open and closed in $H$ and $\alpha_p: Y \ra \bD_p $ is continuous, $j_{\alpha_p}$ defines a homomorphism of $B(\bD_p)$ into $CB(H)$, the algebra of continuous bounded functions on $H$. It is easy to check that $\ja = j_{\alpha_p} \circ j_p$, and $j_p(B(G)) = A_{\omega_p}$, so  $j_{\alpha_p}$ is a homomorphism of $A_{\omega_p}$ into $B(H)$; since $B(H)$ is semisimple, $\ja$ and $j_{\alpha_p}$ are automatically bounded. As noted in the proof of Propostition \ref{CPHomProp}, $j_p^*$ is a normal $*$-isomorphism on $VN_{\omega_p}$, so ${\cal N} = j_p^*(VN_{\omega_p})$ is a von Neumann subalgebra of $W^*(G)$,  and the diagram 
 $$  \xymatrixrowsep{2pc} \xymatrixcolsep{3pc}
\xymatrix{ W^*(H) \ar@<.5ex>[d]_{j_{\alpha_p^*}}  \ar@<.5ex>[rr]^{ \ja^*}
 & & {\cal N} \subseteq W^*(G)
  \\
\ \ \ VN_{\omega_p} \ar@<.5ex>[rru]_{j_p^*} } $$
commutes. The inverse map $(j_p^*)^{-1}$ of $j_p^*: VN_{\omega_p} \ra {\cal N}$ is also a $*$-isomorphism and we have $\ja^* = j_p^* \circ j_{\alpha_p}^*$ and $(j_p^*)^{-1} \circ \ja^* = j_{\alpha_p}^*$, so $\ja^*$ is   positive/ contractive/ completely positive/ completely contractive/ completely bounded   if and only if $j_{\alpha_p}^*$ is so. 
 \end{proof}


 With Corollaries \ref{CPCor1} and \ref{CC Thm Cor}, we  characterized all completely positive and completely contractive  homomorphisms $\vp = \ja : B(G) \ra B(H)$ when $\alpha$ maps into a subgroup $\bD_p$ of $B(G)$.  When  $p$ is further assumed to be a critical idempotent with $\bD_p$ contained in  $G^{\cal E}$, the closure of $G$ in $\DBG$,   we  now observe that we can also sometimes describe when  $\vp$ is positive, contractive, and completely bounded.  Observe that every subgroup of $\DBG$ is abelian---and therefore amenable---when $G$ is abelian and in this case the locally compact subgroups of $G^{\cal E}$ all take the form $\bD_p$ for some critical idempotent $p$, and were characterized by Dunkl and Ramirez in \cite{Dun-Ram}. As we have already stated,  in both  the abelian and nonabelian cases, little is known about homomorphisms of $B(G)$ into $B(H)$.

 \bc  \label{CB Homom for Critical p Cor}  Let $\vp= \ja: B(G) \ra B(H)$ be a  homomorphism such that $\alpha(Y)$ is contained in $\bD_p$ where $p$ is critical,   $\bD_p \subseteq G^{\cal E}$, and  $A_{\omega_p}$ has non-trivial intersection with $A(\bD_p)$. Then:  
 \bi \item[(i)] if $\vp$ is positive/contractive, then $Y$ is an open subgroup/coset in $H$ and the continuous map $\alpha: Y \ra \bD_p \subseteq \DBG$ is either a homomorphism or an anti-homomorphism/affine or anti-affine map; 
 
 \item[(ii)] if $\bD_p$ is amenable and $\vp$ is completely bounded, then $\alpha: Y \subseteq H \ra \DBG$ is a continuous piecewise affine  map.
 \ei  In both cases, $j_{\alpha_p}$ extends to the spatial homomorphism $j_{\alpha_p}: B(\bD_p) \ra B(H)$,  $\vp= \ja$ factors as $\ja = j_{\alpha_p} \circ j_p$ and $j_{\alpha_p}:B(\bD_p) \ra B(H)$ is positive/ contractive/ completely bounded when $\vp$ is so.  \ec

\begin{proof}  As noted above, $A_{\omega_p}$ is a closed, translation-invariant, point-separating, unital subalgebra of $B(\bD_p)$ and---by Proposition \ref{Arsac2.10GenProp}---$A_{\omega_p}$ is conjugation invariant since $\bD_p$ is contained in $G^{\cal E}$. Hence, $A_{\omega_p}$ contains $A(\bD_p)$ by \cite[Corollary 2.3]{Bek-Lau-Sch}.      By Proposition \ref{SpatialHomFactorizationProp}, $\vp = \ja$ factors as $\vp = j_{\alpha_p} \circ j_p$, where $j_{\alpha_p}:A_{\omega_p} \ra B(H)$, $\alpha_p : Y \ra \bD_p$,  is a  spatial homomorphism and  $j_{\alpha_p}$ is positive/contractive in case (i), completely bounded in case (ii).  It follows from Corollary \ref{Lebesgue Decomp Cor}---note that $A_{\omega_p} \supseteq A(\bD_p)$ is required here---that $Y$ and $\alpha_p$ (equivalently $\alpha$) take the desired form in both cases. Moreover,  $j_{\alpha_p}:B(\bD_p) \ra B(H) $ is a positive/contractive homomorphism in case (i) \cite[Corollary 5.6]{Pham} and a completely bounded homomorphism in case (ii) \cite[Corollary 3.2]{Il-Sp}; these statements also follow from our general results in Section 2.    \end{proof}  
 
\br \rm  1. One can similarly obtain  a version of  Corollary \ref{CB Homom for Critical p Cor} for completely positive and completely contractive homomorphisms. However, Corollaries \ref{CPCor1} and \ref{CC Thm Cor} already provide stronger results since they impose no restrictions on $p$. \\
2. One can replace the assumption in Corollary \ref{CB Homom for Critical p Cor} that $\bD_p\subseteq G^{\cal E}$ and $A_{\omega_p}$ has nontrivial intersection with $A(\bD_p)$ with the assumption that $A_{\omega_p} \supseteq A(\bD_p)$. \\
3.  By Corollary \ref{Walter min idempotent etc Corollary}(ii),   $\theta_F$  has dense range in $\bD_{z_F}$ and, as noted on page 125 of \cite{Spr-Sto}, $z_F \in G^{\cal E}$.  It follows that $z_F$ is a critical idempotent such that $\bD_{z_F}$ is contained in $G^{\cal E}$.  Moreover,   since $\bD_{z_F}$ is compact and $A_{\omega_{z_F}}$ is nontrivial, $A_{\omega_{z_F}} = A(\bD_{z_F})$ in this case.  
 \er

\noindent {\bf Acknowledgements:} The author is grateful to Nico Spronk for  offering helpful  points of clarification about the spine compactification. The author is also indebted to Aasaimani Thamizhazhagan for valuable discussions regarding examples satisfying the hypotheses of Corollary \ref{Full characterization of cc and cp in special cases Cor}.

\noindent {\sc Department of Mathematics and Statistics, University
of Winnipeg, 515 Portage Avenue, Winnipeg, MB, R3B 2E9, Canada }

\noindent email: {\tt r.stokke@uwinnipeg.ca}


\begin{thebibliography}{99}


\bibitem{Ars} G. Arsac,  Sur l'espace de Banach engendr$\acute{\rm{e}}$
par les coefficients d'une repr$\acute{\rm{e}}$sentation unitaire,
{\it Publ. D$\acute{e}$p. Math. (Lyon)} 13 (1976), 1-101.

\bibitem{Bek-Lau-Sch} M.E.B. Bekka, A.T.-M. Lau and G.  Schlichting, On invariant subalgebras of the Fourier--Stieltjes algebra of a locally compact group, {\it Math. Ann.} 294 (1992), no. 3, 513-522. 


\bibitem{Ber-Jun-Mil1} J.F. Berglund, H. Junghenn and P.  Milnes, {\it Compact right topological semigroups and generalizations of almost periodicity},  Lecture Notes in Mathematics, 663,  Springer, Berlin, 1978. 

\bibitem{Ber-Jun-Mil} J.F. Berglund, H. Junghenn and P.  Milnes,  {\it Analysis on semigroups: function spaces, compactifications, representations,}  Canadian Mathematical Society Series of Monographs and Advanced Texts,  John Wiley \& Sons, Inc., New York, 1989.




\bibitem{Choi} M.-D. Choi, A Schwarz inequality for positive linear maps on C*-algebras,  {\it Illinois J.  Math.}, 18 (1974), 565-574.

\bibitem{Coh} P. J. Cohen, On homomorphisms of group algebras, {\it
Amer. J. Math} 82 (1960), 213-226.

\bibitem{Dal} H.G. Dales, {\it Banach algebras and automaticcontinuity}, London Math. Soc. Monographs, Volume 24, Clarendon Press, Oxford, (2000).

\bibitem{Dun-Ram} C.F. Dunkl and D.E. Ramirez, Locally compact subgroups of the spectrum of the measure algebra, {\it Semigroup Forum} 3 1971, no. 2, 95-107. 

\bibitem{Eff-Rua} E.G. Effros and Z.-J. Ruan, {\it Operator Spaces}, Oxford
University Press, 2000.




\bibitem{Eym} P. Eymard,   L'alg\'ebre de Fourier d'un groupe localement compact,
 {\it Bull. Soc. Math. France}  92 (1964), 181-236.

\bibitem{FelDor}  J.M.G.  Fell and R.S. Doran, {\it Representations of $*$-algebras,
locally compact groups, and Banach $*$-algebraic bundles  Vol.
1},  Pure and Applied Mathematics, 125. Academic Press, Inc.,
Boston, MA, 1988.

\bibitem{For-Woo} B. Forrest and P. Wood,  Cohomology and the operator space structure of the Fourier algebra and its second dual, {\it Indiana Univ. Math. J.}  50 (2001), no. 3, 1217-1240.

\bibitem{Gli}  I. Glicksberg,  Homomorphisms of certain algebras of measures, {\it Pacific J. Math.} 10 (1960), 167-191.

\bibitem{Gre} F.P. Greenleaf,  Norm decreasing homomorphisms of
 group algebras,  {\it Pacific J. Math.} 15 (1965), 1187-1219.

\bibitem{Hos} B. Host, Le th\'eor\`eme des idempotents dans $B(G)$,
{\it Bull. Soc. Math. France}  114 (1986), no. 2, 215-223.
 
\bibitem{Il} M. Ilie,  On Fourier algebra homomorphisms, {\it J. Funct. Anal.}, 213
(2004), 88-110.


\bibitem{Il-Sp} M. Ilie and N. Spronk,  Completely bounded
homomorphisms of the Fourier algebra, {\it J. Funct. Anal.} 225
(2)(2005), 480-499.



\bibitem{Il-Sp2} M. Ilie and N. Spronk, The spine of a Fourier--Stieltjes algebra, {\it Proc. London Math. Soc.} (3) 94 (2007) 273-301.


\bibitem{Il-St} M. Ilie and R. Stokke,  Weak$^*$-continuous homomorphisms of Fourier-Stieltjes algebras,
{\it  Math. Proc. Cambridge Philos. Soc.},  145 (2008),  107-120.


\bibitem{Kan} E. Kaniuth,  {\it A course in commutative Banach algebras}, Graduate Texts in Mathematics, Springer, New York, 2009.

\bibitem{Kan-Lau} E. Kaniuth and A.T.-M. Lau, \it Fourier and Fourier-Stieltjes algebras on locally compact groups, \rm  Mathematical Surveys and Monographs, 231, American Mathematical Society, Providence, RI, 2018.  

\bibitem{Kan-Lau-Sch} E. Kaniuth, A.T.-M. Lau and G. Schlichting, Lebesgue-type decomposition of subspaces of Fourier--Stieltjes algebras,  {\it Trans. Amer. Math. Soc.} 355 (2003), no. 4, 1467-1490. 

\bibitem{Kan-Lau-Ulg} E. Kaniuth, A.T.-M. Lau and A. $\ddot{\rm{U}}$lger, The Rajchman algebra $B_0(G)$ of a locally compact group $G$, \it  Bull. Sci. Math., \rm 140 (2016), no. 3, 273-302.



\bibitem{Liu-Mis} J.R. Liukkonen and M.W.  Mislove, Fourier--Stieltjes algebras of compact extensions of nilpotent groups, {\it J. Reine Angew. Math.} 325 (1981), 210-220.

\bibitem{Miao} T. Miao, Decomposition of $B(G)$, {\it Trans. Amer. Math. Soc.} 351 (1999), no. 11, 4675-4692. 

\bibitem{Pau} V. Paulsen {\it Completely bounded maps and operator algebras},  Cambridge Studies in Advanced Mathematics, 78, Cambridge University Press, Cambridge, 2002.


\bibitem{Pham} H.L. Pham, Contractive homomorphisms of the
Fourier algebras, {\it Bull. London Math. Soc.} (2010) 42(5), 937-947.




\bibitem{Rud} W. Rudin, {\it Fourier analysis on groups}, Tracts in
Pure and Applied Mathematics, No. 12 Wiley  Interscience,  New
York-London 1962.


\bibitem{Run} V. Runde,  Amenability for dual Banach algebras, {\it Studia Math.} (1) 148 (2001), 47-66.

\bibitem{Run2} V. Runde, Applications of operator spaces to abstract harmonic analysis, {\it Expo. Math.} 22 (2004), no. 4, 317-363.

\bibitem{Run-Spr0} V. Runde and N. Spronk,  Operator amenability of Fourier--Stieltjes algebras,  {\it Math. Proc. Cambridge Philos. Soc.}  136 (2004), no. 3, 675-686.

\bibitem{Run-Spr} V. Runde and N. Spronk,  Operator amenability of Fourier-Stieltjes algebras II, {\it Bull. Lond. Math. Soc.} 39 (2007), no. 2, 194-202.

\bibitem{Spr} N. Spronk, Amenability properties of Fourier algebras and Fourier-Stieltjes algebras: a survey, {\it Banach algebras 2009}, 365-383, Banach Centre Publ., 91, Polish Acad. Sci. Inst. Math., Warsaw, 2010. 


\bibitem{Spr-Sto} N. Spronk and R. Stokke,  Matrix coefficients of unitary representations and associated compactifications of locally compact groups, {\it Indiana Univ. Math. J. }   62 (2013), no. 1, 99-148.
 
   \bibitem{Sto1} R. Stokke, Homomorphisms of convolution algebras, {\it J. Funct. Anal.} 261 (2011), no. 12, 3665-3695.

\bibitem{Sto} R. Stokke, Contractive homomorphisms of measure algebras and Fourier algebras, {\it Studia Math. } 209 (2012), no. 2, 135-150.

\bibitem{Tak} M. Takesaki, \it Theory of Operator Algebras I, \rm Encyclopedia of Mathematical Sciences Vol. 124, Springer--Verlag Berlin Heidelberg, 2002. 

\bibitem{Wal1} M. Walter,  $W\sp{*} $-algebras and nonabelian harmonic
analysis, {\it  J. Funct. Anal.} 11 (1972), 17--38.

\bibitem{Wal2} M. Walter, 
On the structure of the Fourier-Stieltjes algebra, {\it
Pacific J. Math.} 58 (1975), no. 1, 267-281. 

\end{thebibliography}
\end{document}